\documentclass[a4paper,10pt]{article}
\usepackage{amsfonts}
\usepackage{amsmath}
\usepackage{amssymb}
\usepackage{theorem}
\usepackage{url}
\usepackage{amscd}
\usepackage{pstricks}
\usepackage{graphicx}\usepackage{epsfig}
\usepackage[english]{babel}
\input epsf.tex
\newtheorem{thm}{Theorem}[section]
\newtheorem{lem}[thm]{Lemma}

\newtheorem{df}[thm]{Definition}
\newtheorem{ex}[thm]{Example}
\newtheorem{re}[thm]{Remark}
\newtheorem{prop}[thm]{Proposition}

\newcommand{\Irr}[1]{\operatorname{Irr}(#1)}\newcommand{\res}[1]{\operatorname{res}(#1)}

\def\[{[\![}
\def\]{]\!]}
\title{Crystal graphs of higher level $q$-deformed Fock spaces, Lusztig $a$-values and Ariki-Koike algebras}
\author{Nicolas JACON \footnote{Laboratoire de Math\'ematiques de Besan\c{c}on, Universit\'e de Franche-Comt\'e} \footnote{email: jacon@math.univ-fcomte.fr }}
\date{}
\begin{document}
\maketitle
\begin{abstract}
We show that the different labelings of the crystal graph for irreducible
highest weight $\mathcal{U}_q (\widehat{\mathfrak{sl}}_e)$-modules lead to
different labelings of the simple modules for non semisimple Ariki-Koike
algebras by using Lusztig $a$-values.

\end{abstract}

\maketitle

\section{Introduction}
The Ariki-Koike algebras have been introduced by Ariki and Koike in \cite{AK}.
They can be seen as natural generalizations of Iwahori-Hecke algebras of type
$A_{n-1}$ and $B_n$ and as  special cases of cyclotomic Hecke algebras introduced
by  Brou\'e and Malle in \cite{BM}.  Let  $R$ be a commutative ring, let $l\in{\mathbb{N}}$,
$n\in{\mathbb{N}}$ and let  $v,\ x_1,\ x_2,\ ...,\ x_{l}$ be  $l+1$ parameters in $R$.
The Ariki-Koike algebra $\mathcal{H}_{R,n}:=\mathcal{H}_{R,n}(v;x_1,...,x_{l})$
(or cyclotomic Hecke algebra  of type $G(l,1,n)$) over $R$  is the unital associative $R$-algebra presented by:
\begin{itemize}
\item generators: $T_0$, $T_1$,..., $T_{n-1}$,
\item relations: \begin{align*}
               & T_0 T_1 T_0 T_1=T_1 T_0 T_1 T_0,\\
               & T_iT_{i+1}T_i=T_{i+1}T_i T_{i+1}\ (i=1,...,n-2),\\
               & T_i T_j =T_j T_i\ (|j-i|>1),\\
                  &(T_0-x_1)(T_0-x_2)...(T_0-x_{l})  =  0,\\
                             &(T_i-v)(T_i+1)  =  0\ (i=1,...,n-1). \end{align*}
\end{itemize}
Assume that $R$ is a field of characteristic $0$. When $\mathcal{H}_{R,n}$ is a semisimple algebra,
 it is known that the  simple  $\mathcal{H}_{R,n}$-modules are given by the set of
  Specht modules $S^{ \underline{\lambda}}_{R}$ parametrized by the $l$-partitions of rank $n$.
   Using results of Dipper and Mathas \cite{DM}, the study of the
 representation theory of $\mathcal{H}_{R,n}$ in the non semisimple case can be reduced to
   the case where $R=\mathbb{C}$ and :
$$v=\eta_e,\qquad{x_i =\eta_e^{s_j},\ \ j=1,...,l,}$$
where $\displaystyle \eta_e:=\textrm{exp}{(\frac{2i\pi}{e})}\in\mathbb{C}$  and
where $s_j\in{\mathbb{Z}}$ for $j=1,2,...,l$.

In \cite{Ac}, Ariki has first provided a description of the simple modules in
this modular case. Proving and generalizing a conjecture by Lascoux, Leclerc
and Thibon, he has shown in \cite{Ad} that the associated decomposition
matrices can be identified with the Kashiwara-Lusztig canonical basis elements
of the irreducible $\mathcal{U}({\widehat{\mathfrak{sl}}_e})$-module of highest
weight $\Lambda_{s_1(\textrm{mod}\ e)}+\Lambda_{s_2 (\textrm{mod}\
e)}+...+\Lambda_{s_{l}(\textrm{mod}\ e)}$ (where the $\Lambda_i$ denote the
fundamental weights). It is known that these elements can be labeled by using
the ``crystal graph theory'' and there are several ways to do that. Each of
these ways corresponds to a certain realization of the Fock space as a module
over the quantum group $\mathcal{U}_q({\widehat{\mathfrak{sl}}_e})$.  In
particular, one of these realizations  leads to a labeling of the canonical
basis elements by the ``Kleshchev $l$-partitions'' and Ariki has given a
parametrization of the $\mathcal{H}_{\mathbb{C},n}$-simple modules by using
this class of $l$-partitions.

This paper is a continuation of the works done by Geck \cite{G}, Geck and Rouquier \cite{GR}
and  the author \cite{Jp}. The aim is to give natural parametrizations of the
simple modules for non semisimple Ariki-Koike algebras. Let ${\bf
s}_l=(s_1,s_2,...,s_l)\in{\mathbb{Z}^l}$ and let $y$ be an inderterminate. We
consider the Ariki-Koike algebra $\mathcal{H}_{A,n}$  over
$A:=\mathbb{C}[y,y^{-1}]$ with the following choice of parameters:
 \begin{eqnarray*}
&&u_j=y^{lm^{(j)}}\eta_l^{j-1} \textrm{ for}\ j=1,...,l,\\
&&v=y^l.
\end{eqnarray*}
where $\displaystyle \eta_l:=\textrm{exp}(\frac{2i\pi}{l})$ and where for
$j=1,...,l$, we have $\displaystyle
m^{(j)}=s_j-\frac{(j-1)e}{l}+\alpha e$ ($\alpha$ is a positive
integer such that $m^{(j)}\geq 0$ for $j=1,...,l$).

If we specialize $y$ to $\eta_{le}:=\textrm{exp}(\frac{2i\pi}{le})$ via
 a ring homomorphism $\theta$, we obtain the above Ariki-Koike algebra
 $\mathcal{H}_{\mathbb{C},n}:=\mathcal{H}_{\mathbb{C},n}(\eta_e;\eta_e^{s_1},...,\eta_e^{s_l})$. Let $\mathcal{H}_{\mathbb{C}(y),n}:=\mathbb{C}(y)\otimes_A \mathcal{H}_{A,n}$. Then
 $\mathcal{H}_{\mathbb{C}(y),n}$ is a split semisimple algebra and the simple
 $\mathcal{H}_{\mathbb{C}(y),n}$-modules are the Specht modules $S^{\underline{\lambda}}_{\mathbb{C}(y)}$
 defined over $\mathbb{C}(y)$.  We obtain a well-defined decomposition map $d$ between the
 Grothendieck groups of finitely generated $\mathcal{H}_{\mathbb{C}(y),n}$-modules
 and $\mathcal{H}_{\mathbb{C},n}$-modules. For $V\in{\Irr{\mathcal{H}_{\mathbb{C}(y),n}}}$, we have equations:
$$d([V])=\sum_{M\in{\Irr{\mathcal{H}_{\mathbb{C},n}}}}d_{V,M}[M].$$
Using the  symmetric algebra strucure of $\mathcal{H}_{A,n}$, we can attach an integer $a_{{\bf
s}_l}(V)$ to each simple $\mathcal{H}_{\mathbb{C}(y),n}$-module $V$. This is called the
$a$-value of $V$. Note that for Hecke algebras of type $A_{n-1}$ and $B_n$,
this value can also be defined using the Kazhdan-Lusztig theory but such a
theory  is not available
 for the wider case of Ariki-Koike algebras. The aim of this paper is to show
 the following theorem (see Theorem \ref{maing}).\\
\\
{\bf Main Theorem.} {\it For each $M\in{\Irr{\mathcal{H}_{\mathbb{C},n}}}$, there exists a
unique simple $\mathcal{H}_{\mathbb{C}(y),n}$-module $V_M$ such that the following two
conditions hold:
\begin{itemize}
\item $d_{V_M,M}=1$.
\item if there exists $W\in{\Irr{\mathcal{H}_{\mathbb{C}(y),n}}}$ such that $d_{W,M}\neq
  0$ then $a_{{\bf s}_l}(W)>a_{{\bf s}_l}(V_M)$.
\end{itemize}
 The function which sends $M$ to $V_M$ is injective. As a consequence the associated
 decomposition matrix is unitriangular
 and the following set is in natural bijection with $\Irr{\mathcal{H}_{\mathbb{C},n}}$:
$$\mathcal{B}=\{V_M\ |\ M\in{\Irr{\mathcal{H}_{\mathbb{C},n}}}\}.$$
Moreover $\mathcal{B}$ is parametrized by the crystal of the associated
$\mathcal{U}_q(\widehat{\mathfrak{sl}}_e)$-module in one of the realizations of
the Fock space mentioned above.} \vspace{0.3cm}

This result both generalizes and gives a new proof of the main result of
\cite{Jp} where the case $0\leq s_1 \leq ... \leq s_l<e$ was
considered. As the Ariki-Koike algebra $\mathcal{H}_{\mathbb{C},n}$
only depends on the classes modulo $e$ of the integers $s_j$
($j=1,...,l$), the theorem provides several natural
parametrizations for the simple modules of
$\mathcal{H}_{\mathbb{C},n}$. Note also that in the context of Hecke
algebras, the existence of such ``basic sets'' $\mathcal{B}$ is linked with the
existence of a Kazhdan-Lusztig theory (see \cite{Gs}). Hence, this result
gives indices for the existence of such a theory for Ariki-Koike algebras.  

The paper is organized as follows. In section $2$ and $3$, we summarize known
results on the representation theory of Ariki-Koike algebras. We introduce
combinatoric objects and we recall works of Uglov on higher level $q$-deformed
Fock space which will be crucial for the proof of the main result. Section $4$
contains the main theorem of the paper, Theorem \ref{maing}.  The proof of this
Theorem  requires a combinatorial study of the deep results of Uglov
\cite{U}.\\
\\
{\it Acknowledgements.} I thank Bernard Leclerc and Xavier Yvonne for
precious remarks and discussions.

\section{Decomposition maps for Ariki-Koike algebras}\label{dec}

Let  $R$ be a commutative associative ring with unit and let  $v$, $x_1$,..., $x_{l}$ be $l+1$ invertible elements in $R$. Let $n\in{\mathbb{N}}$. Let  $\mathcal{H}_{R,n}:=\mathcal{H}_{R,n}{({v};{x}_1,...,{x}_{l})}$ be the associated Ariki-Koike algebra as it is defined in the introduction. For a complete survey of the representation theory  of $\mathcal{H}_{R,n}$, see \cite{Ma}.

It is known that this algebra is a ``cellular'' algebra in the sense of Graham and Lehrer \cite{GL} and thus has ``Specht modules''  which are parametrized by the $l$-partitions of rank $n$.  A $l$-partition  $\underline{\lambda}$ of rank $n$ is a sequence of $l$ partitions  $\underline{\lambda}=(\lambda^{(1)},...,\lambda^{(l)})$ such that  $\displaystyle{\sum_{k=1}^{l}{|\lambda^{(k)}|}}=n$. We denote by $\Pi^n_l$ the set of  $l$-partitions of rank $n$.

For each $l$-partition $\underline{\lambda}$ of rank $n$, we can associate a remarkable $\mathcal{H}_{R,n}$-module $S_R^{\underline{\lambda}}$ which is free over $R$. This is called  a Specht module (see the definition of ``dual'' Specht modules in \cite{DJM}). Assume that $R$ is a field. In general, the Specht modules are reducible and each of these modules can be endowed with a natural bilinear form. Let $\textrm{rad}(.)$ denotes the radical of this form. For $\underline{\lambda} \in{\Pi_l^n}$, we denote :
$$D^{\underline{\lambda}}_R:= S_R^{\underline{\lambda}}/\textrm{rad}(S_R^{\underline{\lambda}}).$$
Then the theory of cellular algebras gives the following result:
\begin{thm}(Graham-Lehrer \cite{GL}, Dipper-James-Mathas \cite{DJM}) Assume that $R$ is a field then:
\begin{enumerate}
\item Non zero $D_R^{\underline{\lambda}}$ form a complete set of non-isomorphic simple  $\mathcal{H}_{R,n}$-modules. Moreover, all of these modules are absolutely irreducible.
\item If  $\mathcal{H}_{R,n}$ is a semisimple algebra then it is split semisimple and   for all $\underline{\lambda} \in{\Pi_l^n}$ we have  $\textrm{rad}(S_R^{\underline{\lambda}})=0$. Thus, the  $S_R^{\underline{\lambda}}$ form a complete set of non isomorphic simple modules.
\end{enumerate}
\end{thm}

Using this theorem, when  $\mathcal{H}_{R,n}$    is semisimple, the simple
modules are explicitly known and are given by the Specht modules. The following
theorem gives a criterion of semisimplicty for the algebra
$\mathcal{H}_{R,n}$.

\begin{thm}\label{semi} (Ariki \cite{As})\label{semisimple}  $\mathcal{H}_{R,n}$ is split semisimple if and only if:
\begin{itemize}
\item for all $i\neq{j}$ and for all   $d\in{\mathbb{Z}}$ such that $|d|<{n}$, we have:
$$v^d x_i\neq{x_j},$$
\item $\displaystyle{\prod_{i=1}^{n}{(1+v+...+v^{i-1})}\neq{0}}.$
\end{itemize}
\end{thm}

Hence we are reduced to understand the representations of
$\mathcal{H}_{R,n}$ in the modular case.  Assume that $R$ is 
the field of complex numbers $\mathbb{C}$. Then, using the results of Dipper and Mathas (\cite{DM}) and Mathas (\cite{Ms}), the non semisimple case  can be  reduced to the case  where all the $u_i$ are powers of the same root of unity $\eta_e:=\textrm{exp}(\frac{2i\pi}{e})$ with $e\geq 2$. Let  $\mathcal{H}_{\mathbb{C},n}$ be the Ariki-Koike algebra over $\mathbb{C}$ with the following choice of parameters:
\begin{eqnarray*}
&&x_j=\eta_e^{v_{j}}\ \textrm{for}\ j=1,...,l,\\
&&x=\eta_{e},
\end{eqnarray*}
 where $0\leq{v_l}\leq ...\leq v_{1}<e$. The problem of describing the simple modules of  $\mathcal{H}_{\mathbb{C},n}$  is linked with the problem of determining the decomposition map which we now define.

Let $\mathcal{H}_{A,n}$ be an Ariki-Koike defined over a commutative ring $A$ with unit. Let $K$ be the field of fractions of $A$. We assume that:
\begin{itemize}
\item[(1)] $A$ is integrally closed in $K$,
\item[(2)]  $\mathcal{H}_{K,n}:=K\otimes_A \mathcal{H}_{A,n}$ is split semisimple,
\item[(3)] we have a ring homomorphism $\theta: A\to \mathbb{C}$ such that $\mathbb{C}=\textrm{Frac}(\theta (A))$ and such that the specialized algebra $\mathbb{C}\otimes_A \mathcal{H}_{A,n}$ is the Ariki-Koike algebra  $\mathcal{H}_{\mathbb{C},n}$ with the above choice of parameters.
\end{itemize}
Then, by  \cite{Gm}, we have a well-defined decomposition map $d_{\theta}$ between the Grothendieck groups of finitely generated $\mathcal{H}_{K,n}$ and $\mathcal{H}_{\mathbb{C},n}$-modules.
In the context of cellular algebras, this function can be easily defined. Let:
$$\Phi^n_{\{e;v_1,...,v_l\}}:=\{\underline{\lambda} \in{\Pi_l^n}\ |\ D^{\underline{\lambda}}_{\mathbb{C}}\neq 0\}.$$

 Let  $R_0 (\mathcal{H}_{\mathbb{C},n})$  be  the Grothendieck group of finitely generated $\mathcal{H}_{\mathbb{C},n}$-modules. This is generated by the classes of simple $\mathcal{H}_{\mathbb{C},n}$-modules  $[D^{\underline{\lambda}}_{\mathbb{C}}]$ with $\underline{\lambda} \in{\Pi_l^n}$. Hence, for all $\underline{\lambda} \in{\Pi_l^n}$, there exist numbers $d_{\underline{\lambda},\underline{\mu}}$ with $\underline{\mu} \in{\Phi^n_{\{e;v_1,...,v_l\}}}$ such that:
$$[S^{\underline{\lambda}}_{\mathbb{C}}]=\sum_{\underline{\mu} \in{\Phi^n_{\{e;v_1,...,v_l\}}}}{d_{\underline{\lambda},\underline{\mu}}[D_{\mathbb{C}}^{\underline{\mu}}]}.$$
 The matrix $(d_{\underline{\lambda},\underline{\mu}})_{\underline{\lambda}\in{\Pi_l^n},\underline{\mu}\in{\Phi^n}}$ is called the decomposition matrix of $\mathcal{H}_{\mathbb{C},n}$ (where we denote $\Phi^n:=\Phi^n_{\{e;v_1,...,v_l\}}$). Let $\mathcal{F}^n$ be the $\mathbb{C}$-vector space which is generated by the symbols $\[S^{\underline{\lambda}}\]$ with $\underline{\lambda}\in{\Pi_l^n}$. We obtain a homomorphism:
$$\begin{array}{lll}
\mathcal{F}^n&\to&{ R_0 (\mathcal{H}_{\mathbb{C},n})}\\
 \[S^{\underline{\lambda}}\]&\mapsto&[S_{\mathbb{C}}^{\underline{\lambda}}]=\sum_{\underline{\mu}\in{\Phi^n}}{d_{\underline{\lambda},\underline{\mu}}[D_{\mathbb{C}}^{\underline{\mu}}]}.
\end{array}$$

  Now,  $\mathcal{F}^n$ can be naturally identified with the Grothendieck group of finitely generated modules over a semisimple algebra  $\mathcal{H}_{K,n}$ verifying $(1)-(3)$, by identifying the classes of simple  $\mathcal{H}_{K,n}$-modules $[S_K^{\underline{\lambda}}]$ with the symbols  $\[S^{\underline{\lambda}}\]$. Hence,  the decomposition map is defined as follows:
$$\begin{array}{llll}
d:& R_0 (\mathcal{H}_{K,n}) &\to&{ R_0 (\mathcal{H}_{\mathbb{C},n})}\\
 &[S^{\underline{\lambda}}_K]&\mapsto&[S_{\mathbb{C}}^{\underline{\lambda}}]=\sum_{\underline{\mu}\in{\Phi^n}}{d_{\underline{\lambda},\underline{\mu}}[D_{\mathbb{C}}^{\underline{\mu}}]}.
\end{array}$$

By results of Ariki and Uglov, the problem of determining the decomposition matrices for Ariki-Koike algebras can be translated to the problem of computing the canonical basis of $q$-deformed Fock spaces. In the next section, we recall these results.

\section{Canonical basis of higher level $q$-deformed Fock spaces}
The higher $q$-deformed Fock spaces have been introduced in
\cite{JM} (see also \cite{TU}). These spaces which are spanned by the set of
``multipartitions'' can be endowed with a structure of integrable
$\mathcal{U}_q (\widehat{sl}_e)$-module. In  \cite{U},  generalizing
works by  Leclerc and Thibon \cite{LT}, Uglov has given a
construction of  canonical bases for these spaces. In this part, we
review this construction following  \cite{U} (see also \cite{Us} and
\cite{Yt}). This will be used for the proof of the main result of this
paper. Then, we explain the links with the representation theory of Ariki-Koike algebras which are given by Ariki's theorem.

\subsection{$q$-wedge products and $q$-deformed Fock spaces}\label{bij}

 Let $q$ and  $z$ be indeterminates and let $l$ and $e$ be  positive integers. Let $V_e$ be an $e$-dimensional vector space over $\mathbb{Q}(q)$ with basis $v_1^{(e)}$, $v_2^{(e)}$,..., $v_e^{(e)}$.  We put $V_{e,l}:=(V_e\otimes V_l)[z,z^{-1}]$. For $a\in{\{1,...,e\}}$, $b\in{\{1,...,l\}}$ and $m\in{\mathbb{Z}}$, we put $k:=a+e(l-b)-elm$ and $u_k:=v^{(e)}_az^mv_b^{(l)}$. Then $V_{e,l}$ is a  $\mathbb{Q}(q)$ vector space with basis $\{u_k\ |\ k\in{\mathbb{Z}}\}$.

The $q$-wedge square $\bigwedge^2 V_{e,l}$ of  $V_{e,l}$ can be viewed as a $q$-deformation of the exterior square of   $V_{e,l}$. This is a  $\mathbb{Q}(q)$-vector space generated by the monomials $u_{k_1}\wedge u_{k_2}$ with $(k_1,k_2)\in{\mathbb{Z}^2}$. A basis of this space is given by the ``ordered'' monomials that is the monomials $u_{k_1}\wedge u_{k_2}$ such that $k_1> k_2$. Any monomial $u_{k_1}\wedge u_{k_2}$  can be expressed as a linear combination of ordered monomials using the following rules $(R1)$, $(R2)$, $(R3)$ and $(R4)$.

Let $k_1\leq k_2$ and for $i=1,2$, put $k_i=a_i+e(l-b_i)-elm_i$ where
 $a_i\in{\{1,...,e\}}$, $b_i\in{\{1,...,l\}}$ and $m_i\in{\mathbb{Z}}$.
 We define  $\alpha$ and $\beta$ to be the unique integers in $[0,el-1]$ such that $\alpha \equiv(a_2 -a_1 )(\textrm{mod }el) $ and $\beta\equiv(e(b_1-b_2))(\textrm{mod }el)$. Then the relations  $(R1)$, $(R2)$, $(R3)$ and $(R4)$ are given as follows:\\
\begin{itemize}
\item[(R1)] if $\alpha=0$ and  $\beta=0$: $$\displaystyle u_{k_1}\wedge u_{k_2}=-u_{k_2}\wedge u_{k_1};$$
\item[(R2)] if $\alpha\neq 0$ and $\beta=0$:
$$\begin{array}{ll}
u_{k_1}\wedge u_{k_2}=&-\displaystyle q^{-1}u_{k_2}\wedge u_{k_1}+\\
                     &\displaystyle +(q^{-2}-1) \sum_{m\geq 0}q^{-2m} u_{k_2-\alpha-elm}\wedge u_{k_1+\alpha+elm}- \\
                     & \displaystyle-(q^{-2}-1)\sum_{m\geq 1}q^{-2m+1} u_{k_2-elm}\wedge u_{k_1+elm};
\end{array}
$$
\item[(R3)] if $\alpha= 0$ and  $\beta\neq0$:
$$\begin{array}{ll}
u_{k_1}\wedge u_{k_2}=&\displaystyle q u_{k_2}\wedge u_{k_1}+\\
 &\displaystyle +(q^{2}-1)  \sum_{m\geq 0}q^{2m} u_{k_2-\beta-elm}\wedge u_{k_1+\beta+elm}+\\ &\displaystyle  +(q^2 -1)\sum_{m\geq 1}q^{2m-1} u_{k_2-elm}\wedge u_{k_1+elm};\end{array}$$

\item[(R4)]  if $\alpha\neq 0$ and  $\beta\neq 0$:
$$\begin{array}{ll}
  u_{k_1}\wedge u_{k_2}=&u_{k_2}\wedge u_{k_1}+\\
&+\displaystyle{(q-q^{-1})\sum_{m\geq{0}} \frac{q^{2m+1}+q^{-2m-1}}{q+q^{-1}}  u_{k_2-\beta-elm}\wedge u_{k_1+\beta+elm}}+\\
&+\displaystyle{(q-q^{-1})\sum_{m\geq{0}} \frac{q^{2m+1}+q^{-2m-1}}{q+q^{-1}}  u_{k_2-\alpha-elm}\wedge u_{k_1+\alpha+elm}}+\\
&+\displaystyle{(q-q^{-1})\sum_{m\geq{0}} \frac{q^{2m}-q^{-2m}}{q+q^{-1}}  u_{k_2-\beta-\alpha-elm}\wedge u_{k_1+\beta+\alpha+elm}+}\\
&+\displaystyle{(q-q^{-1})\sum_{m\geq{1}} \frac{q^{2m}-q^{-2m}}{q+q^{-1}}  u_{k_2-elm}\wedge u_{k_1+elm}};\end{array}$$
\end{itemize}
where the summations continue as long as the monomials appearing under
the sums are ordered.

For any integer $r\geq 2$, we can now define the $r$-fold $q$-wedge product
$\bigwedge^r V_{e,l}$. This is the  $\mathbb{Q}(q)$-vector space generated by
the elements $u_{k_1}\wedge u_{k_2}\wedge ... \wedge u_{k_r}$ with
$k_i\in{\mathbb{Z}}$. Again the ordered monomials, that is the monomials
$u_{k_1}\wedge u_{k_2}\wedge ... \wedge u_{k_r}$  with $k_1>k_2>...>k_r$, form
a basis of $\bigwedge^r V_{e,l}$. Moreover, an arbitrary monomial  can be
expressed  as a linear combination of ordered monomials using the relations
$(R1)$,  $(R2)$,  $(R3)$ and $(R4)$ in every adjacent pair of the factors.

Finally, for  $s\in{\mathbb{Z}}$, the semi-infinite $q$-wedge product $\bigwedge^{s+\frac{\infty}{2}} V_{e,l}$ of charge $s$ is the inductive limit of $\bigwedge^r V_{e,l}$ where the maps $\bigwedge^r V_{e,l}\to\bigwedge^{r+1} V_{e,l}$ are given by $v\mapsto v\wedge u_{s-r}$. Hence,  $\bigwedge^{s+\frac{\infty}{2}} V_{e,l}$ is the  $\mathbb{Q}(q)$-vector space generated by the semi-infinite monomials  $u_{k_1}\wedge u_{k_2}\wedge ...$ with  $k_i\in{\mathbb{Z}}$ and such that $k_i=s-i+1$ for $i>>0$. A basis is given by the ordered monomials  that is the monomials $u_{k_1}\wedge u_{k_2}\wedge ...$ with $k_1>k_2>....$ and the relations  $(R1)$-$(R4)$ provide a way to express an arbitrary monomial as a linear combinaison of ordered monomials.

Now, each semi-infinite ordered monomial can be labeled by a pair
$(\underline{\lambda},{\bf s}_l)$ where $\underline{\lambda}\in{\Pi_l^n}$ is a
$l$-partition of rank $n$ and ${\bf s}_l=(s_1,...,s_l)\in{\mathbb{Z}}^l$ is
such that $\sum_{i=1}^l s_i =s$. This labeling coincides with that of \cite{U}
up to transformation $(s_1,...,s_l)\mapsto (s_l,...,s_1)$ and
$(\lambda^{(1)},...,\lambda^{(l)})\mapsto (\lambda^{(l)},...,\lambda^{(1)})$.
Let  $u_{k_1}\wedge u_{k_2}\wedge ...\in{ \bigwedge^{s+\frac{\infty}{2}}
V_{e,l}}$ be a  semi-infinite ordered monomial. For $i=1,2,...$, we put
$k_i=a_i+e(l-b_i)-elm_i$ where   $a_i\in{\{1,...,e\}}$, $b_i\in{\{1,...,l\}}$
and $m_i\in{\mathbb{Z}}$. For $b=1,....,l$, let $k_1^{(b)}>k_2^{(b)}>...$ be
the semi-infinite sequence obtained by ordering the elements of the set
$\{a_i-em_i\ |\ b_i=b\}$ in strictly decreasing order. Then, there is a unique
$s_b\in{\mathbb{Z}}$ such that $k_i^{(b)}=s_b-i+1$ for $i>>1$. We put
$\underline{\lambda}=(\lambda^{(1)},...,\lambda^{(l)})$ where for $i>0$,
$\lambda^{(b)}_i=k_i^{(b)}-s_b+i-1$ and ${\bf s}_l=(s_1,...,s_l)$. This defines
a bijection between the set of semi-infinite ordered monomials and the set of
pairs  $(\underline{\lambda},{\bf s}_l)$ where
$\underline{\lambda}\in{\Pi_l^n}$ and ${\bf
s}_l=(s_1,...,s_l)\in{\mathbb{Z}}^l$ is such that $\sum_{i=1}^l s_i =s$.

 This bijection  can be described as follows (see \cite{Y} for much
 details).  First, the infinite decreasing sequence $\underline{k}=(k_1,k_2,...)$
 can be pictured as a set of black  beads on an infinite runner. For
 example, the following sequence :
 $$\underline{k}=(15,12,8,7,3,1,-2,-4,-5,-6,...)$$
  is represented by the following abacus:

\begin{center}
\leavevmode
\epsfxsize= 12cm
\epsffile{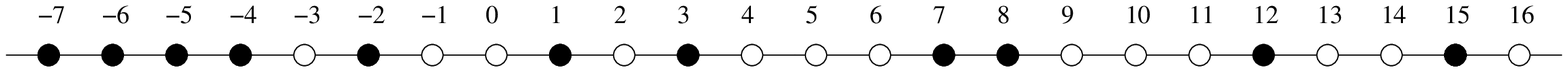}
\end{center}

On the other hand,  using the decomposition $k_i=a_i+e(l-b_i)-elm_i$, the same sequence can be represented by an $l$-abacus, that is, as a set of black  beads on $l$ infinite runners. For $e=4$ and $l=3$, we obtain:

\begin{center}
\leavevmode
\epsfxsize= 12cm
\epsffile{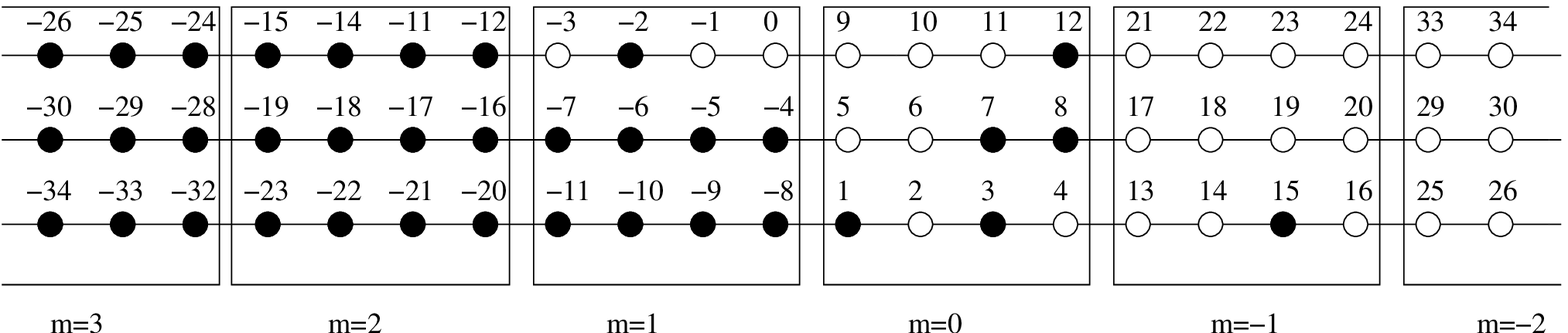}
\end{center}

In the above representation, the beads are labeled by the integers $k_i$ (with $i=1,2,...$). We can alternatively labeled them by the integers $k_i^{(b)}$ (with $i=1,2,...$ and $b=1,2,...,l$) as follows:

\begin{center}
\leavevmode
\epsfxsize= 12cm
\epsffile{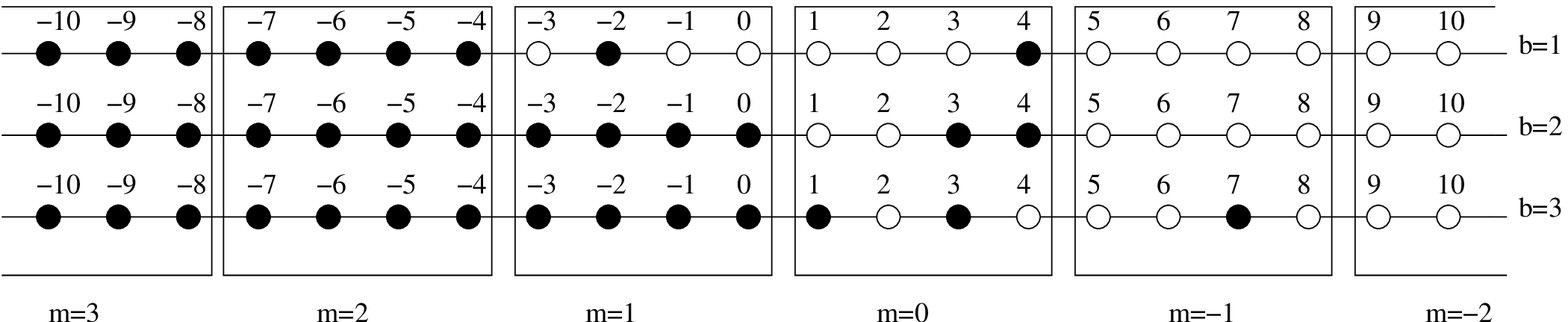}
\end{center}

Let  $\underline{\lambda}$ be the $l$-partition  and ${\bf s}_l=(s_1,...,s_l)$ be the $l$-tuple of integers  associated to $\underline{k}$. For $i=1,2,...$ and $b=1,2,...,l$, we have $\lambda^{(b)}_i=k_i^{(b)}-s_b+i-1$. Thus, one can now easily determine  $\underline{\lambda}$   by counting the number of white beads at the left of each black bead on each runner. Continuing the above example, we obtain $\underline{\lambda}=((6,1),(2,2),(4,1))$ and ${\bf s}_l=(-2,2,3)$.

Now, the higher level $q$-deformed Fock space of charge ${\bf s}_l $ is defined to be the $\mathbb{C}(q)$-vector space generated by the symbols $|\underline{\lambda},{\bf s}_l \rangle$ with $\underline{\lambda}\in{\Pi_l^n}$ :
$$\mathcal{F}_{e,{\bf s}_l}:=\bigoplus_{\overset{\underline{\lambda}\in{\Pi_l^n}}{n\in{\mathbb{N}}}} \mathbb{C}(q)|\underline{\lambda},{\bf s}_l\rangle.$$
Thus, if we identify the semi-infinite ordered monomial $u_{\underline{k}}\in{\bigwedge^{s+\frac{\infty}{2}} V_{e,l}}$  with the pair  $|\underline{\lambda},{\bf s}_l\rangle$ defined by the above bijection, we get:
$$\bigwedge^{s+\frac{\infty}{2}} V_{e,l}=\bigoplus_{s_1+...+s_l=s}\mathcal{F}_{e,{\bf s}_l}.$$

\subsection{Canonical basis and Ariki's theorem}\label{ariki}

Let $\mathcal{U}_q({\widehat{\mathfrak{sl}}}_e)$ be the quantum group
 associated to  the Lie algebra ${\widehat{\mathfrak{sl}}}_e$. We
 denote by $e_i$, $f_i$ and $k_i$ with $i=0,...,e-1$, the Chevalley
 generators and by  $\mathcal{U}_q'({\widehat{\mathfrak{sl}}}_e)$ the
 subalgebra of $\mathcal{U}_q({\widehat{\mathfrak{sl}}}_e)$ generated
 by $e_i$, $f_i$, $k_i$, $k_i^{-1}$. We can construct a structure of
 $\mathcal{U}_q'({\widehat{\mathfrak{sl}}}_e)$-module on
 $\bigwedge^{s+\frac{\infty}{2}}V_{e,l}$ (see \cite[\S 2.1,\S 3.5, \S 4.2]{U})
such that each of the subspace $\mathcal{F}_{e,{\bf s}_l}$ is stable with respect to this action.
Before describing the action of the  $\mathcal{U}_q'({\widehat{\mathfrak{sl}}}_e)$-module $\mathcal{F}_{e,{\bf s}_l}$,  we need some combinatorial definitions.

  Let   $\underline\lambda={(\lambda^{(1)} ,...,\lambda^{(l)})}$ be a $l$-partition of rank  $n$. The diagram of  $\underline{\lambda}$ is the following set:
$$[\underline{\lambda}]=\left\{ (a,b,c)\ |\ 1\leq{c}\leq{l},\ 1\leq{b}\leq{\lambda_a^{(c)}}\right\}.$$
The elements of this diagram are called   the nodes of  $\underline{\lambda}$.
Let  $\gamma=(a,b,c)$ be a node of  $\underline{\lambda}$. The residue
of  $\gamma$  associated to the set  $\{e;{s_1},...,{s_{l}}\}$ is the
unique element of $[0,e[$ such that:
$$\textrm{res}{(\gamma)}\equiv (b-a+s_{c})(\textrm{mod}\ e).$$
If $\gamma $ is a node with residue $i$, we say that  $\gamma$ is an  $i$-node.
Let  $\underline{\lambda}$ and  $\underline{\mu}$ be two $l$-partitions of rank  $n$ and $n+1$ such that  $[\underline{\lambda}]\subset{[\underline{\mu}]}$. There exists a node $\gamma$ such that  $[\underline{\mu}]=[\underline{\lambda}]\cup{\{\gamma\}}$. Then, we denote $[\underline{\mu}]/[\underline{\lambda}]=\gamma$ and  if $\textrm{res}{(\gamma)}=i$, we say that  $\gamma$ is an addable $i$-node  for   $\underline{\lambda}$ and a removable  $i$-node for  $\underline{\mu}$. Now, we introduce an order on the set of nodes of a $l$-partition.
 We say that  $\gamma=(a,b,c)$ is {\it above}  $\gamma'=(a',b',c')$ if:
$$b-a+s_c<b'-a'+s_{c'}\ \textrm{or } \textrm{if}\ b-a+s_c=b'-a'+s_{c'}\textrm{ and }c'<c.$$
Let $\underline{\lambda}$ and $\underline{\mu}$ be two $l$-partitions of rank $n$ and $n+1$ such that  there exists an $i$-node $\gamma$ such that  $[\underline{\mu}]=[\underline{\lambda}]\cup{\{\gamma\}}$. We define the following numbers:

\begin{align*}
{N}_i^{a}{(\underline{\lambda},\underline{\mu})}=&   \sharp\{ \textrm{addable }\ i-\textrm{nodes of } \underline{\lambda}\ \textrm{ above } \gamma\} \\
 & -\sharp\{ \textrm{removable }\ i-\textrm{nodes of } \underline{\mu}\ \textrm{ above } \gamma\},\\
{N}_i^{b}{(\underline{\lambda},\underline{\mu})}= &   \sharp\{ \textrm{addable } i-\textrm{nodes of } \underline{\lambda}\ \textrm{ below } \gamma\}\\
&   -\sharp\{ \textrm{removable } i-\textrm{nodes of } \underline{\mu}\ \textrm{ below } \gamma\},\\
{N}_{i}{(\underline{\lambda})} =& \sharp\{ \textrm{addable } i-\textrm{nodes of } \underline{\lambda}\}\\
& -\sharp\{ \textrm{removable } i-\textrm{nodes of } \underline{\lambda}\}.
\end{align*}
Now, $\mathcal{F}_{e,{\bf s}_l}$ is an integrable $\mathcal{U}_q'(\widehat{\mathfrak{sl}_e})$-module with action:
\begin{align*}
e_{i}|\underline{\lambda},{\bf s}_l\rangle &=\sum_{\res{[\underline{\lambda}]/[\underline{\mu}]}\equiv i}{q^{-{N}_i^{a}{(\underline{\mu},\underline{\lambda})}}}|\underline{\mu},{\bf s}_l\rangle,\\
f_{i}|\underline{\lambda},{\bf s}_l\rangle &=\sum_{\res{[\underline{\mu}]/[\underline{\lambda}]}\equiv i}{q^{{N}_i^{b}{(\underline{\lambda},\underline{\mu})}}}|\underline{\mu},{\bf s}_l\rangle,\\
k_{i}|\underline{\lambda},{\bf s}_l\rangle&=q^{{N}_i{(\underline{\lambda})}}|\underline{\lambda},{\bf s}_l\rangle.
\end{align*}
where $0\leq{i}\leq{e-1}$. 
Note that this action coincides with that of \cite{U} up to transformation $(s_1,...,s_l)\mapsto (s_l,...,s_1)$.

We now introduce another basis for $\mathcal{F}_{e,{\bf s}_l}$ namely  the Leclerc-Thibon canonical basis. This basis  is defined by using an involution on the wedge space  $\bigwedge^{s+\frac{\infty}{2}} V_{e,l}$ which has been introduced in \cite{LT} for $l=1$ and generalized to any $l$ in \cite{U}.
 Let $u_{\underline{k}}:=u_{k_1}\wedge u_{k_2}\wedge ...$ be a semi infinite monomial (ordered or not). For $i=1,2,...$, we put
 $k_i=a_i+e(l-b_i)-elm_i$ where $a_i\in{\{1,...,e\}}$,
 $b_i\in{\{1,...,l\}}$ and $m_i\in{\mathbb{Z}}$. For $r\in{\mathbb{N}}$, we put:
$$\omega(u_{\underline{k}})=\sharp \{i<j\ |\ a_i=a_j\},$$
$$\omega'(u_{\underline{k}})=\sharp \{i<j\ |\ b_i=b_j\}.$$
Then, for $r\geq \sum_{i=1}^{\infty} k_i-(s-i+1)$ we set:
$$\overline{u_{\underline{k}}}:=(-q)^{\omega'(u_{\underline{k}})}q^{-\omega(u_{\underline{k}})} u_{k_r}\wedge u_{k_{r-1}}\wedge .... \wedge u_{k_1}\wedge u_{k_{r+1}}\wedge u_{k_{r+2}}\wedge ...$$
One can prove that  this monomial is independent of $r$ and that
$u\mapsto{\overline{u}}$ defines an involution on
$\bigwedge^{s+\frac{\infty}{2}} V_{e,l}$. The canonical basis is now defined as
follows:
\begin{thm}(Leclerc-Thibon \cite[Theorem  4.1]{LT}, Uglov \cite[\S 4.4]{U})\label{cb} Let $s\in{\mathbb{Z}}$. There exists a unique basis:
 $$\{G(\underline{\lambda},{\bf s}_l)\ |\ \sum_{i=1}^l s_i=s,\ \underline{\lambda}\in{\Pi_l^n},\ n\in{\mathbb{N}}\}$$
of $\bigwedge^{s+\frac{\infty}{2}} V_{e,l}$ such that:
\begin{itemize}
\item $\overline{ G(\underline{\lambda},{\bf s}_l)}=G(\underline{\lambda},{\bf s}_l)$,
\item $G(\underline{\lambda},{\bf s}_l) - |\underline{\lambda},{\bf s}_l\rangle\in \bigoplus_{\underline{\mu}} q\mathbb{C}[q] |\underline{\mu},{\bf s}_l\rangle   $.
\end{itemize}
This is called the canonical basis of $\bigwedge^{s+\frac{\infty}{2}} V_{e,l}$.
 \end{thm}
In particular the set $\{G(\underline{\lambda},{\bf s}_l)\ |\ \underline{\lambda}\in{\Pi_l^n},\ n\in{\mathbb{N}}\}$ gives a basis
 of the Fock space $\mathcal{F}_{e,{\bf s}_l}$. Now, we consider the subspace $M_{{\bf s}_l}:=\mathcal{U}_q'(\widehat{\mathfrak{sl}_e}).|\underline{\emptyset},{\bf s}_l\rangle$.
It is well known that this is isomorphic to the irreducible $\mathcal{U}_q'(\widehat{\mathfrak{sl}_e})$-module $V(\Lambda)$ with highest weight
$\Lambda:=\Lambda_{s_{1}(\textrm{mod}\ e)}+\Lambda_{s_{2}(\textrm{mod}\ e)}+...+\Lambda_{s_{l}(\textrm{mod}\ e)}$.
Note that if ${\bf s}_l'=(s_1',...,s_l')\in{\mathbb{Z}^l}$ is such that $s_i\equiv s_i'(\textrm{mod }e)$, then the modules
 $M_{{\bf s}_l}$ and $M_{{\bf s}_l'}$ are isomorphic (but the action of  $\mathcal{U}_q'(\widehat{\mathfrak{sl}_e})$ on the elements
 of the standard basis $|\underline{\lambda},{\bf s}_l\rangle$ and  $|\underline{\lambda},{\bf s}_l'\rangle$ are different in general).

A basis of $M_{{\bf s}_l}$ can be given by using the canonical basis of $\mathcal{F}_{e,{\bf s}_l}$
 and by studying the associated crystal graph. This graph can be described combinatorially as follows.

Let   $\underline{\lambda}$ be a   $l$-partition and  let $\gamma$ be  an $i$-node of $\underline{\lambda}$,
 we say that  $\gamma$ is a normal $i$-node of  $\underline{\lambda}$ if,
whenever $\eta$ is an $i$-node of $\underline{\lambda}$   below  $\gamma$,
there are more removable $i$-nodes between $\eta$ and $\gamma$ than addable $i$-nodes between  $\eta$ and $\gamma$.
If $\gamma$ is the highest normal $i$-node of    $\underline{\lambda}$, we say that  $\gamma$ is a good $i$-node. Note that this notion depends on the choice of ${\bf s}_l$.

Then, the crystal graph of $\mathcal{F}_{e,{\bf s}_l}$ is given by:
\begin{itemize}
\item vertices: the $l$-partitions,
\item edges: $\displaystyle{{\underline{\lambda}\overset{i}{\rightarrow}{{\underline{\mu}}}}}$ if and only if $[\underline{\mu}]/[\underline{\lambda}]$ is a good  $i$-node.
\end{itemize}

By using properties of crystal bases, we can obtain the crystal graph of
$M_{{\bf s}_l}$: this is the connected components of that of
$\mathcal{F}_{e,{\bf s}_l}$ which contain the vacuum vector
$|\underline{\emptyset},{\bf s}_l\rangle$. The vertices of this graph are given
by the following class of $l$-partitions.

\begin{df} Let ${\bf s}_l\in{\mathbb{Z}^n}$.  The set of Uglov $l$-partitions $\Lambda^n_{e;{\bf s}_l}$ is defined   recursively  as follows.
\begin{itemize}
   \item We have $\underline{\emptyset}:=(\emptyset,\emptyset,...,\emptyset)\in{  \Lambda^n_{e;{\bf s}_l} }$.
    \item If $\underline{\lambda}\in\Lambda^n_{e;{\bf s}_l}$ , there exist $i\in{\{0,...,e-1\}}$ and a good $i$-node $\gamma$ such that if we remove $\gamma$ from  $\underline{\lambda}$, the resulting  $l$-partition is in $\Lambda^n_{e;{\bf s}_l}$.
\end{itemize}
\end{df}

\begin{re}
Assume that  ${\bf s}_l\in{\mathbb{Z}^n}$ is such that $0\leq s_1 \leq
s_2\leq ...\leq s_{l}<e$ then it is shown in \cite[Proposition 2.11]{FL} that the set of Uglov $l$-partitions are the  $l$-partitions  $\underline\lambda={(\lambda^{(1)} ,...,\lambda^{(l)})}$ such that:
\begin{enumerate}
\item for all $1\leq{j}\leq{l-1}$ and $i=1,2,...$, we have:
\begin{align*}
&\lambda_i^{(j)}\geq{\lambda^{(j+1)}_{i+s_{j+1}-s_j}},\\
&\lambda^{(l)}_i\geq{\lambda^{(1)}_{i+e+s_1-s_{l}}};
\end{align*}

\item  for all  $k>0$, among the residues appearing at the right ends of the length $k$ rows of   $\underline\lambda$, at least one element of  $\{0,1,...,e-1\}$ does not occur.
\end{enumerate}
 Such $l$-partitions are called FLOTW $l$-partitions in \cite{Jp}.

Assume that   ${\bf s}_l\in{\mathbb{Z}^n}$ is such that $s_1>> s_2>>
...>> s_{l}$, then the set  Uglov $l$-partitions $\Lambda^n_{e;{\bf
    s}_l}$ coincides with the set of  Kleshchev $l$-partitions by
\cite[\S 3]{Ac}.
\end{re}
Now, the canonical basis of $M_{{\bf s}_l}$ is the following set:
$$\{G(\underline{\lambda},{\bf s}_l)\ |\ \underline{\lambda}\in{\Lambda^{n}_{e;{\bf s}_l}},\ n\in{\mathbb{N}}\}.$$

The following theorem gives a link between the canonical basis elements of  $M_{{\bf s}_l}$ and the decomposition matrices of Ariki-Koike algebras.

\begin{thm}\label{Ar} (Ariki, see \cite[Theorem 4.49]{Ab})
  Let  $\mathcal{H}_{\mathbb{C},n}$ be the Ariki-Koike algebra over $\mathbb{C}$ with the following choice of parameters:
\begin{eqnarray*}
&&x_j=\eta_e^{v_{j}}\ \textrm{for}\ j=1,...,l,\\
&&v=\eta_{e},
\end{eqnarray*}
where $0\leq v_1\leq ...\leq v_l<e$. Let $\Phi^n:=\Phi^n_{\{e;v_l,...,v_l\}}$ and   $(d_{\underline{\lambda},\underline{\mu}})_{\underline{\lambda}\in{\Pi_l^n},\underline{\mu}\in{\Phi^n}}$ be the associated decomposition matrix.

Let  ${\bf s}_l=(s_1,...,s_l)\in{\mathbb{Z}^l}$ be such that $v_j\equiv s_j(\textrm{mod}\ e)$ for $j=1,...,l$. Then, for each $\underline{\nu}\in{\Lambda^n_{ e;{\bf s}_l   }}$, write
$$G(\underline{\nu},{\bf s}_l)=\sum_{\underline{\lambda}\in{\Pi^n_l}}d_{\underline{\lambda},\underline{\nu}}(q)|\underline{\lambda},{\bf s}_l\rangle,$$
 where $d_{\underline{\lambda},\underline{\nu}}(q)\in{\mathbb{C}[q]}$.
 Then there exists a  bijection $\kappa:{\Phi^n}\to \Lambda^n_{e;{\bf s}_l}$ such that for all $\underline{\lambda}\in{\Pi_l^n}$ and  $\underline{\mu}\in{\Phi^n}$, we have:
$$d_{\underline{\lambda},\underline{\mu}}=d_{\underline{\lambda},\kappa(\underline{\mu})}(1).$$
\end{thm}

In other words, Ariki's theorem asserts that the columns of the decomposition
matrix of  $\mathcal{H}_{\mathbb{C},n}$ with the above choice of parameters
coincides with the canonical basis elements of  $M_{{\bf s}_l}$ evaluated at
$q=1$ whenever $v_j\equiv s_j(\textrm{mod}\ e)$ for $j=1,...,l$. Note that this
decomposition matrix only depends on $\{e;v_1,...,v_l\}$ whereas  the canonical
basis elements of  $M_{{\bf s}_l}$ depends on $\{e;s_1,...,s_l\}$. In
particular, if $s_j\equiv s_j'(\textrm{mod}\ e)$ for  $j=1,...,l$,  the
labelings of the canonical basis elements  of  $M_{{\bf s}_l}$ and $M_{{\bf
s}_l'}$ by   $\Lambda^n_{e; {\bf s}_l}$ and $\Lambda^n_{e; {\bf s}_l'}$ are
different in general as we can see in the following example.
\begin{ex}\label{exem}
Assume that $l=2$, $e=4$, $v_1=0$, $v_2=1$. Then different values for ${\bf
s}_l$ lead to different labelings of the same crystal graph which are
given in Example \ref{exfin}.
\end{ex}

 In \cite{Ac}, Ariki has given an explicit description of $\Phi^n_{\{e;v_1,...,v_l\}}$ by showing that this set coincides with the set   $\Lambda^n_{ e;{\bf s}_l}$ with  $ s_1>> s_2>> ...>> s_{l}$ and  $v_j\equiv s_j(\textrm{mod}\ e)$ for $j=1,...,l$. In \cite{Jp}, another parametrization of the simple   $\mathcal{H}_{\mathbb{C},n}$-modules has been given by using the set  $\Lambda^n_{e; {\bf s}_l}$ with ${\bf s}_l=(v_1,...,v_l)$ (namely the set of FLOTW $l$-partitions) and an ordering of the rows of the decomposition matrices by Lusztig $a$-values. In the next section, we will show that each of the sets  $\Lambda^n_{ e;{\bf s}_l}$ with    $v_j\equiv s_j(\textrm{mod}\ e)$ has a natural interpretation in  the representation theory of  $\mathcal{H}_{\mathbb{C},n}$.

\section{Unitriangularity of the decomposition matrices of Ariki-Koike algebras}
 \subsection{Specialisations and Lusztig $a$-values}
Let $e$ be an integer and let ${\bf s}_l=(s_1,s_2,...,s_{l})\in{\mathbb{Z}^l}$ be a sequence of integers.  The aim of this part is to study the Ariki-Koike algebra $\mathcal{H}_{\mathbb{C},n}$ with the following choice of parameters:
\begin{eqnarray*}
&&u_j=\eta_e^{s_j(\textrm{mod }e)}\ \textrm{for}\ j=1,...,l,\\
&&v=\eta_e.
\end{eqnarray*}
For $j=1,...,l$, we define rational numbers:
$$m^{(j)}=s_j-\frac{(j-1)e}{l}+\alpha e,$$
where $\alpha$ is a positive  integer such that $m^{(j)}\geq{0}$  for $j=1,...,l$. Let $y$ be an indeterminate and put $A:=\mathbb{C}[y,y^{-1}]$. We consider the Ariki-Koike algebra $\mathcal{H}^{{\bf s}_l}_{A,n}$ over $A$   with the following parameters:
\begin{eqnarray*}
&&u_j=y^{lm^{(j)}}\eta_l^{j-1}\ \textrm{for}\ j=1,...,l,\\
&&v=y^l
\end{eqnarray*}
 (this algebra is  denoted by
$\mathcal{H}_{A,n}$ in the introduction but here, we use the
notation  $\mathcal{H}_{A,n}^{{\bf s}_l}$ to insist on the
dependance on ${\bf s}_l$).

By Theorem \ref{semi},   $\mathcal{H}_{\mathbb{C}(y),n}^{{\bf s}_l}:=\mathbb{C}(y)\otimes_A \mathcal{H}_{A,n}^{{\bf s}_l}$ is split semisimple.
Moreover, if we specialize the parameter $y$ to ${\eta_{le}:=\textrm{exp}(\frac{2i\pi}{le})}$,
we obtain the above Ariki-Koike algebra  $\mathcal{H}_{\mathbb{C},n}$. Hence, we have a well-defined
decomposition map between $R_0 (\mathcal{H}_{\mathbb{C}(y),n}^{{\bf s}_l})$  and  $R_0 (\mathcal{H}_{\mathbb{C},n})$.

We now associate to each simple $\mathcal{H}_{\mathbb{C}(y),n}^{{\bf s}_l}$-module $S^{\underline{\lambda}}_{\mathbb{C}(y)}$ an $a$-value following \cite{Jp}. Put ${\bf m}_l=(m^{(1)},...,m^{(l)})$ where the $m^{(j)}$ are defined above. We need to define the
notion of ``${\bf m}_l$-translated symbols'' associated to $l$-compositions.

Let $n\in{\mathbb{N}}$ and $l\in{\mathbb{N}}$. An $l$-composition $\underline{\lambda}$ of rank $n$  is an $l$-tuple $(\lambda^{(1)},....,\lambda^{(l)})$ where :
\begin{itemize}
\item for all $i=1,...,l$, we have $\lambda^{(i)}=(\lambda^{(i)}_1,...,\lambda^{(i)}_{h^{(i)}})$ for  $h^{(i)}\in{\mathbb{N}}$ and $\lambda^{(i)}_j\in{\mathbb{N}_{>0}}$ ($j=1,...,h^{(i)}$), $h^{(i)}$ is called the height  of $\lambda^{(i)}$,
\item $\displaystyle{\sum_{i=1}^{l}\sum_{j=1}^{h^{(i)}}\lambda^{(i)}_j=n}$.

\end{itemize}
 Let $\underline{\lambda}=(\lambda^{(1)},...,\lambda^{(l)})$ be a $l$-composition and let $h^{(i)}$ be the heights of the compositions $\lambda^{(i)}$. Then the height of $\underline{\lambda}$ is the following positive integer:
$$h_{\underline{\lambda}}:=\textrm{max}{\{h^{(1)},...,h^{(l)}\}}.$$
Let $k$ be a positive integer. The translated symbol $\textbf{B}[{\bf m}_l]'$ associated to ${\bf m}_l$, $k$ and $\underline{\lambda}$  is:
$$\textbf{B}[{\bf m}_l]':=(B'^{(1)},...,B'^{(l)}),$$
where $B'^{(i)}$, for $i=1,...,l$, is given by:
$$B'^{(i)}:=(B'^{(i)}_1,...,B'^{(i)}_{h_{\underline{\lambda}}+k}),$$
in which:
$$B'^{(i)}_{j}:=\lambda^{(i)}_j-j+h_{\underline{\lambda}}+k+m^{(i)}\  \  (1\leq j \leq h_{\underline{\lambda}}+k).$$
The integer  $h_{\underline{\lambda}}+k$ is called the height of $\textbf{B}[{\bf m}_l]'$.

Now, the $a$-value on the $\mathcal{H}_{\mathbb{C}(y),n}^{{\bf
    s}_l}$-module $S^{\underline{\lambda}}_{\mathbb{C}(y)}$ is defined
to be the lower degree of the associated the Schur element. These
    elements have been explicitely comptuted by 
Geck, Iancu and Malle in \cite{GIM} (see \cite[Proposition 3.1]{Jp}). 
We obtain the following definition:

\begin{df} Let  $\underline{\lambda}$ be a $l$-partition of rank $n$  and let  $S_{\mathbb{C}(y)}^{\underline{\lambda}}$ be the simple $\mathcal{H}^{{\bf s}_l}_{\mathbb{C}(y),n}$-module.
Let  $k$ be a positive integer and $\textbf{B}[{\bf m}_l]'$ be the
${\bf m}_l$-translated symbol associated to ${\bf m}_l$ and $k$. Then, if $h$ is the height of  $\textbf{B}[{\bf m}_l]'$, the $a$-value of  $S_{\mathbb{C}(y)}^{\underline{\lambda}}$ is the following rational number:
$$a_{{\bf s}_l}(\underline{\lambda})=f(n,h,{\bf m}_l)+\sum_{{1\leq{i}\leq{j}<l+1}\atop{{(a,b)\in{B'^{(i)}\times{B'^{(j)}}}}\atop{a>b\ \textrm{if}\ i=j}}}{\min{\{a,b\}}}   -\sum_{{1\leq{i,j}<l+1}\atop{{a\in{B'^{(i)}}}\atop{1\leq{k}\leq{a}}}}{\min{\{k,m^{(j)}\}}},$$
where $f(n,h,{\bf m}_l)$ is a rational number which only depends on
the parameters $\{e;s_1,...,s_{l}\}$, on $h$ and on $n$ (the
expression of $f$ is given in \cite[Proposition 3.2]{Jp}).
\end{df}
We need to introduce the following preorder on the set of $l$-compositions.\begin{df}
Let  $\underline{\mu}$ and  $\underline{\nu}$ be  $l$-compositions
of rank  $n$. Let $k$ and $k'$ be such that $h_{\underline{\mu}}+k=h_{\underline{\nu}}+k'$. Let     $\textbf{B}_{ \underline{\mu}}[{\bf m}_l]'$,
(resp.  $\textbf{B}_{ \underline{\nu}}[{\bf m}_l]'$) be the ${\bf
  m}_l$-translated symbol associated  to $\underline{\mu}$, ${\bf m}_l$ and $k$ (resp.   $\underline{\mu}$, ${\bf m}_l$ and $k'$). Then we write:
$$\underline{\mu} \prec{\underline{\nu}},$$
if:
$$\sum_{{1\leq{i}\leq{j}<l+1}\atop{{(a,b)\in{B_{\underline{\mu}}'^{(i)}
\times{B_{\underline{\mu}}'^{(j)}}}}\atop{a>b\ \textrm{if}\ i=j}}}
{\min{\{a,b\}}}   -\sum_{{1\leq{i,j}<l+1}
\atop{{a\in{B_{\underline{\mu}}'^{(i)}}}
\atop{1\leq{k}\leq{a}}}}{\min{\{k,m^{(j)}\}}}<$$
$$\sum_{{1\leq{i}\leq{j}<l+1}\atop{{(a,b)\in{B_{\underline{\nu}}'^{(i)}
\times{B_{\underline{\nu}}'^{(j)}}}}\atop{a>b\ \textrm{if}\ i=j}}}
{\min{\{a,b\}}}   -\sum_{{1\leq{i,j}<l+1}\atop{{a\in{B_{\underline{\nu}}
'^{(i)}}}\atop{1\leq{k}\leq{a}}}}{\min{\{k,m^{(j)}\}}}.$$
In particular, if  $\underline{\mu}$ and   $\underline{\nu}$ are
 $l$-partitions, as the translated symbol of $\underline{\mu}$ and
$\underline{\nu}$ have the same height $h$, we have:
$$\underline{\mu} \prec{\underline{\nu}}\iff{a_{{\bf s}_l}(\underline{\mu})
<a_{{\bf s}_l}(\underline{\nu})}.$$
\end{df}

 Before beginning the proof of the main result, we need to give  the
 following useful proposition which has been shown in \cite[Lemma 4.4]{Jp} in the case  $0\leq s_1 \leq ....\leq s_l <e$. One can easily check that the proof holds in the general case.

\begin{prop}\label{combi} Let  $\underline{\lambda}$ be a
 $l$-composition of rank $n$, let  $\textbf{B}:=(B^{(1)},...,B^{(l)})$
be an ordinary symbol of $\underline{\lambda}$,  let  $\beta_1$
and $\beta_2$ be two elements of   $\textbf{B}[{\bf m}_l]'$, we assume
 that:
$$\beta_1 < \beta_2.$$
 Let  $r\in{\mathbb{N}}$. We add  $r$ nodes  to $\underline{\lambda}$
 on the part associated to    $\beta_1$. Let  $\underline{\mu}$
 be the resulting $l$-composition of rank $n+r$.  We add  $r$ nodes
 to $\underline{\lambda}$  on the part associated to    $\beta_2$.
 Let  $\underline{\nu}$ be the resulting $l$-composition of rank
$n+r$.   Then, we have:
$$\underline{\nu} \prec{\underline{\mu}}.$$
In particular, if   $\underline{\mu}$ and   $\underline{\nu}$
are some  $l$-partitions, as the translated symbols of
$\underline{\mu}$ and   $\underline{\nu}$ have the same height, we have:
$$a_{{\bf s}_l}(\underline{\mu})>a_{{\bf s}_l}(\underline{\nu}).$$
\end{prop}
Assume that  $\underline{\lambda}$ is a  $l$-composition of rank $n$.

\begin{center}
\leavevmode
\epsfxsize= 12cm
\epsffile{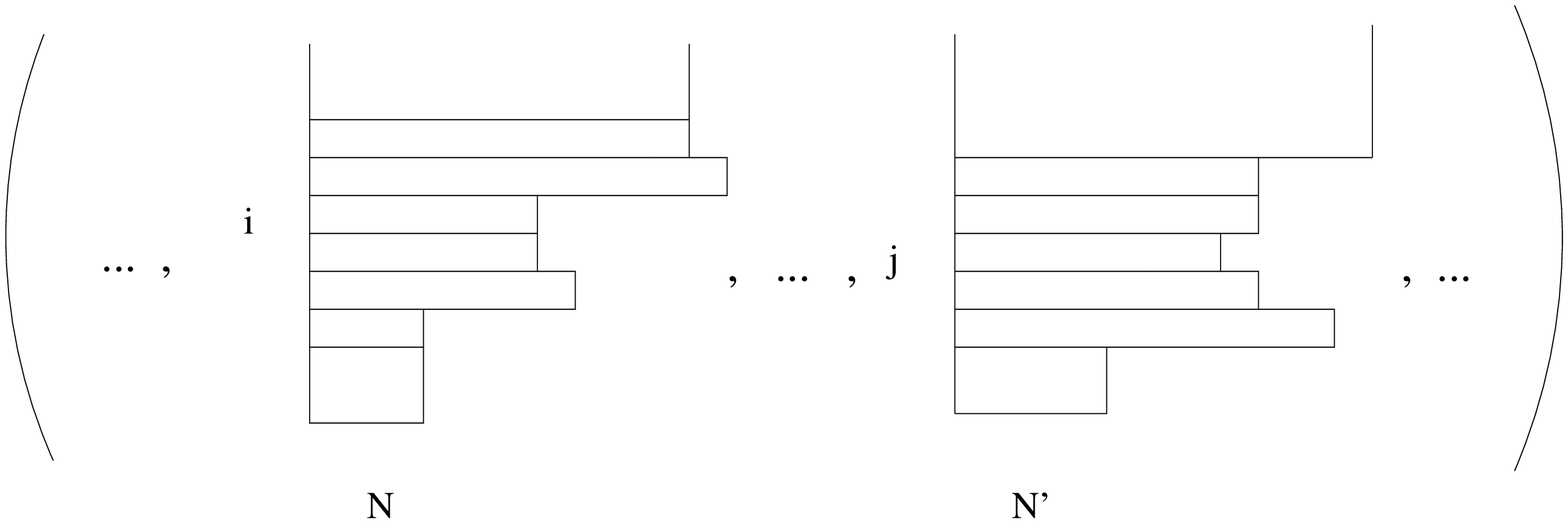}
\end{center}

Assume  that   $\underline{\mu}$ and   $\underline{\mu}'$ are two
$l$-compositions which are obtained from  $\underline{\lambda}$ by
adding $r$ nodes on the parts  $\lambda^{(N)}_i$ and
$\lambda^{(N')}_j$ respectively where $1\leq N, N'\leq l$. Assume in addition that we have:
$$\lambda^{(N)}_i-i+s_N-\frac{(N-1)e}{l}>\lambda^{(N')}_j-j+s_{N'}-\frac{(N'-1)e}{l}.$$
Keeping the above representation of $\underline{\lambda}$, $\underline{\mu}$ and   $\underline{\mu}'$ are respectively given as follows:

\begin{center}
\leavevmode
\epsfxsize= 12cm
\epsffile{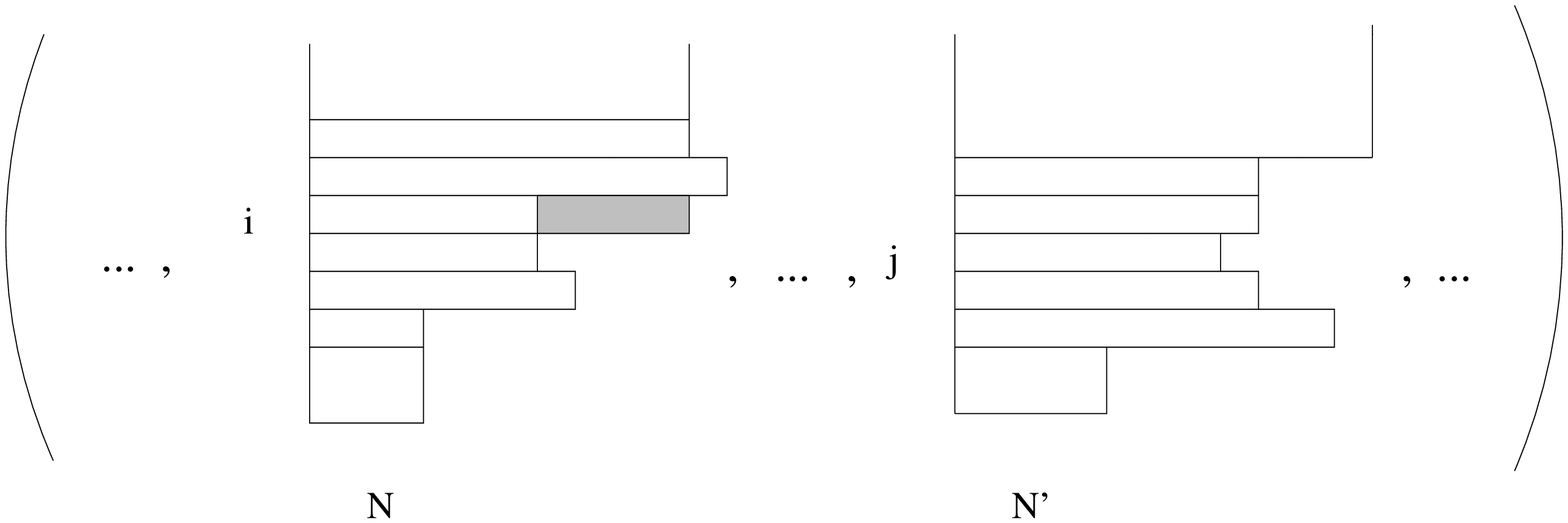}
\end{center}

\begin{center}
\leavevmode
\epsfxsize= 12cm
\epsffile{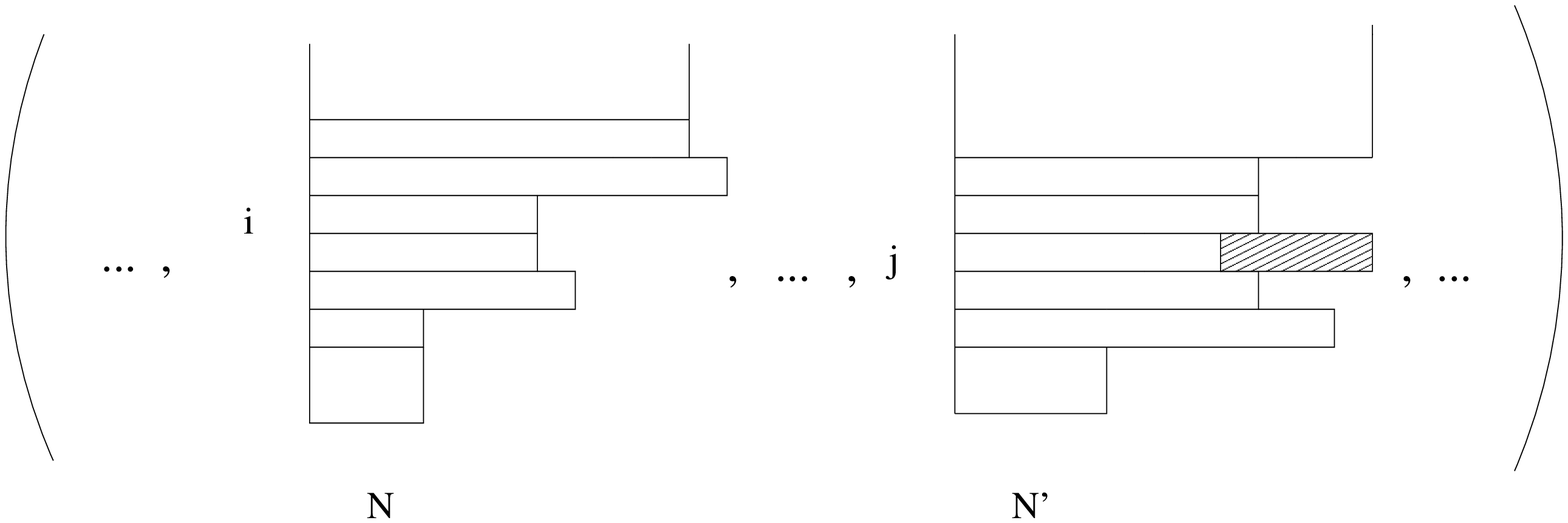}
\end{center}

Then Proposition \ref{combi} asserts that we have:
$$\underline{\mu} \prec{\underline{\mu}'}.$$

The aim of the following part is to show that the matrix associated to the  canonical basis elements of the irreducible modules $M_{{\bf s}_l}$ is unitriangular with respect to $a$-values.

\subsection{Study of the involution on the semi-infinite $q$-wedge product}

Let $s\in{\mathbb{Z}}$. We keep the notations of the previous sections. We will work with the  semi-infinite $q$-wedge product $\bigwedge^{s+\frac{\infty}{2}} V_{e,l}$. First, we attach to each semi-infinite monomial an $a$-value.  Let $u_{\underline{k}}\in \bigwedge^{s+\frac{\infty}{2}} V_{e,l}$ be a   semi-infinite $q$-wedge product and let  $u_{\underline{\tilde{k}}}$ be the monomial obtained by reordering the $k_i$ in strictly decreasing order. Then, the bijection described in section \ref{bij} shows that we can associate  to  $u_{\underline{\tilde{k}}}$ a symbol $|\underline{\lambda},{\bf s}_l\rangle$ with $\underline{\lambda}\in{\Pi_l^n}$ and ${\bf s}_l=(s_1,...,s_l)\in{\mathbb{Z}^l}$ verifying $\sum_{i=1}^l s_i =s$. We put:
\begin{align*}
\pi (u_{\underline{k}})&:=|\underline{\lambda},{\bf s}_l\rangle,\\
a(u_{\underline{k}})&:=a_{{\bf s}_l}(\underline{\lambda}).
\end{align*}
In this part, we show the following proposition:
\begin{prop}\label{main2} Let  $u_{\underline{k}}\in \bigwedge^{s+\frac{\infty}{2}} V_{e,l}$ be a   semi-infinite ordered $q$-wedge product. Then, we have:
$$\overline{u_{\underline{k}}}= u_{\underline{k}}+\textrm{ sum of ordered monomials }u_{\underline{r}}\textrm{ with }a(u_{\underline{r}})>a(u_{\underline{k}}), $$
where the involution $\overline{ }$ is defined in section \ref{ariki}.
\end{prop}
We note that a similar property is shown in \cite{U} but for a  (partial) order which is not compatible with the preorder induced by the $a$-values.

Let  $u_{\underline{k}}$ be an arbitrary semi-infinite monomial such that there exists  $i\in{\mathbb{N}}$ such that $k_i<k_{i+1}$. Then the relations $(R1)-(R4)$ in section \ref{bij} show how to express   $u_{\underline{k}}$ in terms of semi-infinite monomials   $u_{\underline{k}'}$ such that $k_i'>k_{i+1}'$. Let  $u_{\underline{k}'}$ be such a   semi-infinite monomial such that $k_i'\neq k_{i+1}$ and $k_{i+1}'\neq k_i$. Then the  $l$-partitions associated to $u_{\underline{k}}$  and  $u_{\underline{k}'}$  have the same charge ${\bf s}_l$. We put :
\begin{align*}
|\underline{\lambda},{\bf s}_l\rangle&:=\pi (u_{\underline{k}}),\\
|\underline{\lambda}',{\bf s}_l\rangle&:=\pi (u_{\underline{k}'}).
\end{align*}
From the relations $(R1)-(R4)$, using the representations of   $u_{\underline{k}}$  and  $u_{\underline{k}'}$ by $l$-abacus, we can see that  $u_{\underline{k}'}$ is obtained from  $u_{\underline{k}}$ by moving two black beads representing $k_i$ and $k_{i+1}$  to two black beads representing  $k_i'$ and $k_{i+1}'$ and lying in the same runners as $k_i$ and $k_{i+1}$. As a consequence, the $l$-partition $\underline{\lambda}'$ is obtained from $\underline{\lambda}$ by removing a ribbon $R$ of size $x$ from a component $b$ and adding a ribbon $R'$ of size $x$ to a component $b'$ as in the following figure.

\begin{center}
\leavevmode
\epsfxsize= 12cm
\epsffile{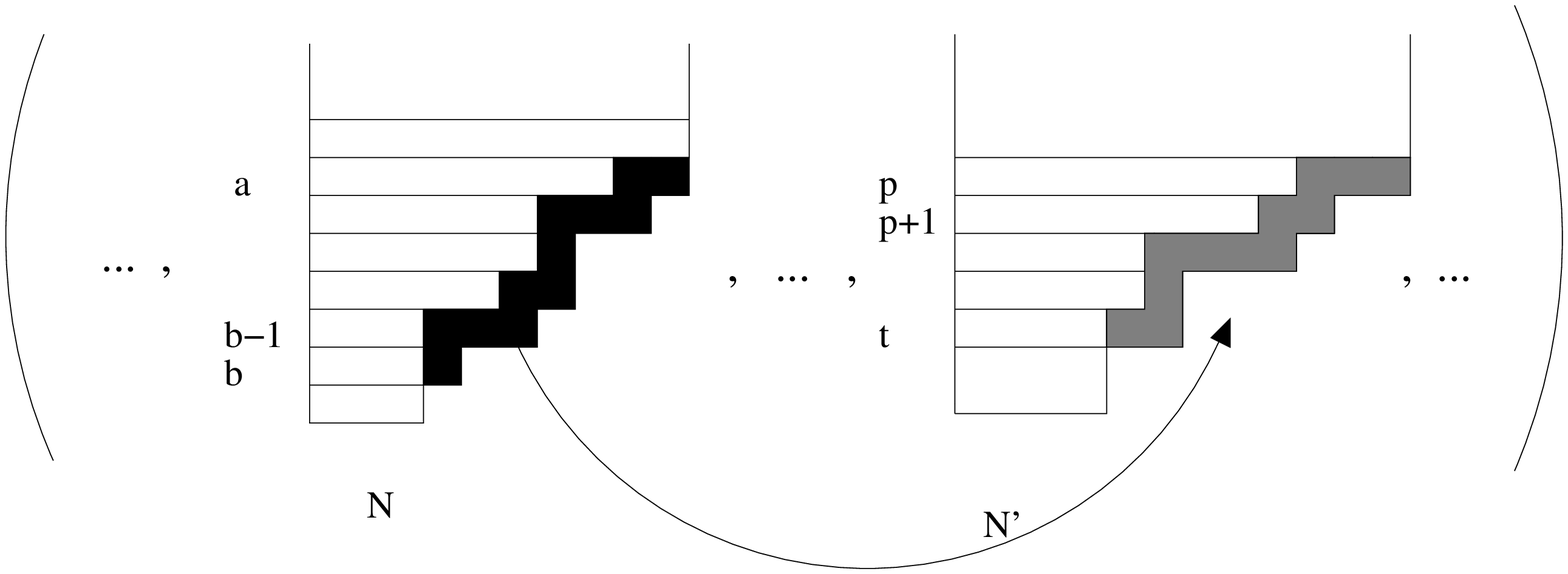}
\end{center}

 Hence, there exist positive integers  $r$, $s$, $p$ and $t$ such that:
$$\lambda'^{(b)}=(\lambda^{(b)}_1,....,\lambda^{(b)}_{r-1},\lambda^{(b)}_{r+1}-1,...,,\lambda^{(b)}_{s}-1,\lambda^{(b)}_r-(x-(r+s)),\lambda^{(b)}_{s+1},...,\lambda_{h^{(b)}}^{(s)}),$$
$$\lambda'^{(b')}=(\lambda^{(b')}_1,...,\lambda_{p-1}^{(b')},\lambda_t^{(b')}+(x-(p-k)),\lambda^{(b')}_{p}+1,...,\lambda^{(b')}_{t-1}+1,\lambda^{(b')}_{t+1},...,\lambda^{(b')}_{h^{(b')}}).   $$
In fact, the two $l$-partitions $\underline{\lambda}'$ and
$\underline{\lambda}$ are both obtained by adding $R$ and $R'$ to the
components $b$ and $b'$ of  the same $l$-partition
$\underline{\nu}=\pi(u_{\underline{r}})$ where
$$\underline{r}=(k_1,...,k_{i-1},\hat{k}_i,\hat{k}_{i+1},k_{i+2},...).$$ 
Here,  $\hat{k}_i$ and $\hat{k}_{i+1}$ are the minimal integers of
$\{k_i,k_{i+1},k_i ',k_{i+1} '\}$ such that the black beads associated to  $\hat{k}_i$ and $\hat{k}_{i+1}$ lie at the components $b$ and $b'$ of the $l$-abacus of  $u_{\underline{r}}$ and   $\hat{k}_{i+1}<\hat{k}_{i}$.

\begin{lem}\label{pcomb} We keep the above notations. 
Assume that $\underline{\lambda}\neq \underline{\lambda}'$. Put
$\hat{k}_i=a+e(l-b)-el m$ and  $\hat{k}_{i+1}=a'+e(l-b')-el m'$ where
$a,a'\in{\{1,...,e\}}$, $b,b'\in{\{1,...,l\}}$ and,
$m,m'\in{\mathbb{Z}}$.    Then we have:
$$a-em-\frac{(b-1)e}{l}>a'-em'-\frac{(b'-1)e}{l}.$$
\end{lem}
\emph{Proof}. For $j=i,\ i+1$, put $k_j=a_j+e(l-b_j)-elm_j$ where  $a_j\in{\{1,...,e\}}$,
$b_j\in{\{1,...,l\}}$ and $m_j\in{\mathbb{Z}}$. As in section
\ref{bij}, we define  $\alpha$
and $\beta$ to be the unique integers in $[0,el-1]$ such that $\alpha
  \equiv(a_{i+1} -a_i )(\textrm{mod }el) $ and
  $\beta\equiv e(b_i-b_{i+1})(\textrm{mod }el)$. As
  $\underline{\lambda}\neq \underline{\lambda}'$, we can assume that
  $\alpha\neq 0$ or $\beta\neq 0$. 

Now, regarding the relations $(R1)-(R4)$ in section \ref{bij}, we have
$5$ cases to consider.

\begin{itemize}
\item  Assume that $k_{i}'=k_{i+1}-elm$ and
  $k_{i+1}'=k_{i}+elm$ with $m> 0$. If $\beta=0$, we have  $\hat{k}_i=k_{i+1}'$
  and $\hat{k}_{i+1}=k_{i }$ and since $k_{i+1}'>k_i$, the result is
  clear. If $\beta\neq 0$,  we have  $\hat{k}_i=k_{i}'$
  and $\hat{k}_{i+1}=k_{i }$.

As $k_{i}'>k_{i+1}'$, we have:
$$a_{i+1}+e(l-b_{i+1})-el(m_{i+1}+m)>a_{i}+e(l-b_{i})-el(m_{i}-m).$$
Now, as $m>0$ and $\displaystyle a_{i+1}-a_i>\frac{a_{i+1}-a_i}{l}-e$, we conclude
that:
$$a_{i+1}-e(m_{i+1}+m)-(b_{i+1}-1)\frac{e}{l}>a_{i}-em_{i}-(b_{i}-1)\frac{e}{l}.$$

\item Assume that $k_{i}'=k_{i+1}-\beta-elm$ and
  $k_{i+1}'=k_{i}+\beta+elm$ with $m\geq 0$ and $\beta\neq 0$. Then, we have
  $\hat{k}_i=k_{i+1}'$    and $\hat{k}_{i+1}=k_{i } $.

\begin{itemize}
 \item If $b_{i+1}<b_{i}$ then $\beta=e(b_i-b_{i+1})$ and we have:
\begin{align*}
\hat{k}_i&=a_{i}+e(l-b_{i})-elm_i+e(b_i-b_{i+1})+elm\\
&=a_i+e(l-b_{i+1})-el(m_i-m).\end{align*}
We conclude that:
$$a_i-e(m_i-m)-(b_{i+1}-1)\frac{e}{l}>a_i-em_i-(b_i-1)\frac{e}{l}.$$
 \item If $b_i<b_{i+1}$  then $\beta=e(b_i-b_{i+1})+el$ and we have:
\begin{align*}
\hat{k}_i&=a_{i}+e(l-b_{i})-elm_i+e(b_i-b_{i+1})+elm+el\\
&=a_i+e(l-b_{i+1})-el(m_i-m-1).\end{align*}
 We conclude that:
$$a_i-e(m_i-m-1)-(b_{i+1}-1)\frac{e}{l}>a_i-em_i-(b_i-1)\frac{e}{l}.$$
\end{itemize}

\item Assume that $k_{i+1}'=k_{i}+\alpha+elm$ and
  $k_{i}'=k_{i+1}-\alpha-elm$ with $m\geq 0$ and $\alpha\neq 0$. If
  $\beta=0$, we have $\hat{k}_i=k_{i+1}'$    and $\hat{k}_{i+1}=k_{i
  }$ and since $k_{i+1}'>k_i$, the result is clear. If $\beta\neq 0$, we have
  $\hat{k}_i=k_{i}'$    and $\hat{k}_{i+1}=k_{i }$.

\begin{itemize}
 \item If $a_{i+1}-a_i>0$ then $\alpha=a_{i+1}-a_i$ and we have:
\begin{align*}
\hat{k}_i&=a_{i+1}+e(l-b_{i+1})-elm_{i+1}-a_{i+1}+a_{i}+elm\\
&=a_i+e(l-b_{i+1})-el(m_{i+1}-m).\end{align*}
We conclude that:
$$a_i-e(m_{i+1}-m)-(b_{i+1}-1)\frac{e}{l}>a_i-em_{i}-(b_i-1)\frac{e}{l}.$$
Indeed, if $m_{i+1}-m<m_i$, the result is clear. If otherwise, as
$m_{i+1}\leq m_i$, we have $m=0$ and $m_{i+1}=m_i$. Hence, as
$k_{i}'>k_{i+1}'$, we have $b_{i+1}<b_i$.
\item  If $a_{i+1}-a_i<0$ then $\alpha=a_{i+1}-a_i+el$ and we have:
\begin{align*}
\hat{k}_i&=a_{i+1}+e(l-b_{i+1})-elm_{i+1}-a_{i+1}+a_{i}+el(m+1)\\
&=a_i+e(l-b_{i+1})-el(m_{i+1}-m-1).\end{align*}
We conclude that:
$$a_i-e(m_{i+1}-m-1)-(b_{i+1}-1)\frac{e}{l}>a_i-em_{i}-(b_i-1) \frac{e}{l}.$$
\end{itemize}
\item Assume that $k_{i+1}'=k_{i}+\alpha+\beta+elm$ and
  $k_{i}'=k_{i+1}-\alpha-\beta-elm$ with $m\geq 0$, $\beta\neq 0$ and
  $\alpha\neq 0$. Then, we have
  $\hat{k}_i=k_{i+1}'$    and $\hat{k}_{i+1}=k_{i }$.
\begin{itemize}
\item if $a_{i+1}>a_i$ and $b_{i}>b_{i+1}$ then $\alpha=a_{i+1}-a_i$
  and $\beta=e(b_{i}-b_{i+1})$. We have:
\begin{align*}
\hat{k}_i&=a_{i}+e(l-b_{i})-elm_i+e(b_i-b_{i+1})+(a_{i+1}-a_i)+elm\\
&=a_{i+1}+e(l-b_{i+1})-el(m_i-m).\end{align*}
We conclude that:
$$a_{i+1}-e(m_i-m)-(b_{i+1}-1)\frac{e}{l}>a_i-em_i-(b_i-1)\frac{e}{l}.$$
\item if $a_{i+1}<a_i$ and $b_{i}>b_{i+1}$ then $\alpha=a_{i+1}-a_i+el$
  and $\beta=e(b_{i}-b_{i+1})$. We have:
\begin{align*}
\hat{k}_i&=a_{i}+e(l-b_{i})-elm_i+e(b_i-b_{i+1})+(a_{i+1}-a_i)+el(m+1)\\
&=a_{i+1}+e(l-b_{i+1})-el(m_i-m-1).\end{align*}
We conclude that:
$$a_{i+1}-e(m_i-m-1)-(b_{i+1}-1)\frac{e}{l}>a_i-em_i-(b_i-1)\frac{e}{l}.$$
\item if $a_{i+1}>a_i$ and $b_{i}<b_{i+1}$ then $\alpha=a_{i+1}-a_i$
  and $\beta=e(b_{i}-b_{i+1})+el$. We have:
\begin{align*}
\hat{k}_i&=a_{i}+e(l-b_{i})-elm_i+e(b_i-b_{i+1})+(a_{i+1}-a_i)+elm+el\\
&=a_{i+1}+e(l-b_{i+1})-el(m_i-m-1).\end{align*}
We conclude that:
$$a_{i+1}-e(m_i-m-1)-(b_{i+1}-1)\frac{e}{l}>a_i-em_i-(b_i-1)\frac{e}{l}.$$
\item if $a_{i+1}>a_i$ and $b_{i}<b_{i+1}$ then $\alpha=a_{i+1}-a_i+el$
  and $\beta=e(b_{i}-b_{i+1})+el$. We have:
\begin{align*}
\hat{k}_i&=a_{i}+e(l-b_{i})-elm_i+e(b_i-b_{i+1})+(a_{i+1}-a_i)+elm+2el\\
&=a_{i+1}+e(l-b_{i+1})-el(m_i-m-2).\end{align*}
We conclude that:
$$a_{i+1}-e(m_i-m-2)-(b_{i+1}-1)\frac{e}{l}>a_i-em_i-(b_i-1)\frac{e}{l}.$$
\end{itemize} 
  \end{itemize}
\rightline{$\Box$}

We are now able to prove Proposition \ref{main2}. Keeping the notations
of the beginning of this section,  Lemma \ref{pcomb} shows that we have:
$$\nu_{s}^{(b)}-s+s_{b}-\frac{(b-1)e}{l}>\nu_{t}^{(b')}-t+s_{b'}-\frac{(b'-1)e}{l}.\ \ \ \ \ \ (*)$$

Assume that we have $\lambda_i'^{(b')}=\nu_i^{(b')}+y_i$ for integers $y_i$ with $i=p,p+1,...,t$ and assume that we have   $\lambda_j^{(b)}=\nu_j^{(b)}+x_j$ for integers $x_j$ with $j=r,r+1,...s$.
We define a $l$-composition $\underline{\lambda}[1]$ as follows.
\begin{itemize}
\item If $y_t\leq x_s$, then we put $j_1:=s$ and we define:
$$\underline{\lambda}[1]_i^{(k)}=\left\{\begin{array}{ll}
{\lambda}_i'^{(k)}&\textrm {if }(i,k)\neq (s,b) \textrm{ and }(i,k)\neq (t,b'),\\
{\lambda}_i'^{(k)}+y_t&\textrm {if }(i,k)= (s,b),\\
{\lambda}_i'^{(k)}-y_t&\textrm {if }(i,k)= (t,b').\\
\end{array}\right.$$
In this case, $\underline{\lambda}[1]$ and $\underline{\lambda}'$ are both obtained by adding the same number of nodes on a $l$-composition. Then, by $(*)$ and by Proposition \ref{combi}, we get:
$$\underline{\lambda}[1]\prec\underline{\lambda}'.$$

\item If  $y_t> x_s$, there exists $j_1\in{\{r,r+1,...,s-1\}}$ such that:
$$x_{j_1}+x_{j_1+1}+....+x_s\geq y_t \geq x_{j_1+1}+...+x_{s-1}+x_s.$$
Then, we define:
$$ \underline{\lambda}[1]_i^{(k)}=\left\{  \begin{array}{ll}
{\lambda}_i'^{(k)}+x_i&\textrm {if }k=b\textrm{ and }\\
                      & \textrm{   }s\geq i\geq j_1+1,\\
{\lambda}_i'^{(k)}+y_t-(x_{j_1+1}+...+x_{s-1}+x_s)&\textrm {if }(i,k)=(j_1,b),\\
{\lambda}_i'^{(k)}-y_t&\textrm {if }(i,k)= (t,b'),\\
{\lambda}_i'^{(k)}&\textrm {if otherwise}.\\
\end{array}\right.$$
For all $i$ such that $r\leq i\leq s$, observe that
$\nu^{(b)}_s+x_{i+1}+...+x_{s-1}+x_s=\nu^{(b)}_{i}+(s-i)$. Then, by $(*)$, for all $i$ such that $j_1\leq i\leq s$, we have:
$$\nu_{i}^{(b)}-i+s_{b}-\frac{(b-1)e}{l}>\nu_{t}^{(b')}-t+x_{i+1}+...+x_{s}+s_{b'}-\frac{(b'-1)e}{l}.$$
Thus,  for all $i$ such that $j_1\leq i\leq s$, we obtain:
$${\lambda'}_{i}^{(b)}-i+s_{b}-\frac{(b-1)e}{l}>{\lambda'}_{t}^{(b')}-t+x_{i+1}+...+x_{s-1}+x_s-y_t+s_{b'}-\frac{(b'-1)e}{l}.$$

Thus, using Proposition \ref{combi}, we get:
$$\underline{\lambda}[1]\prec\underline{\lambda}'.$$
\end{itemize}
Keeping the above figure, $\underline{\lambda}[1]$ is  as follows.
\begin{center}
\leavevmode
\epsfxsize= 12cm
\epsffile{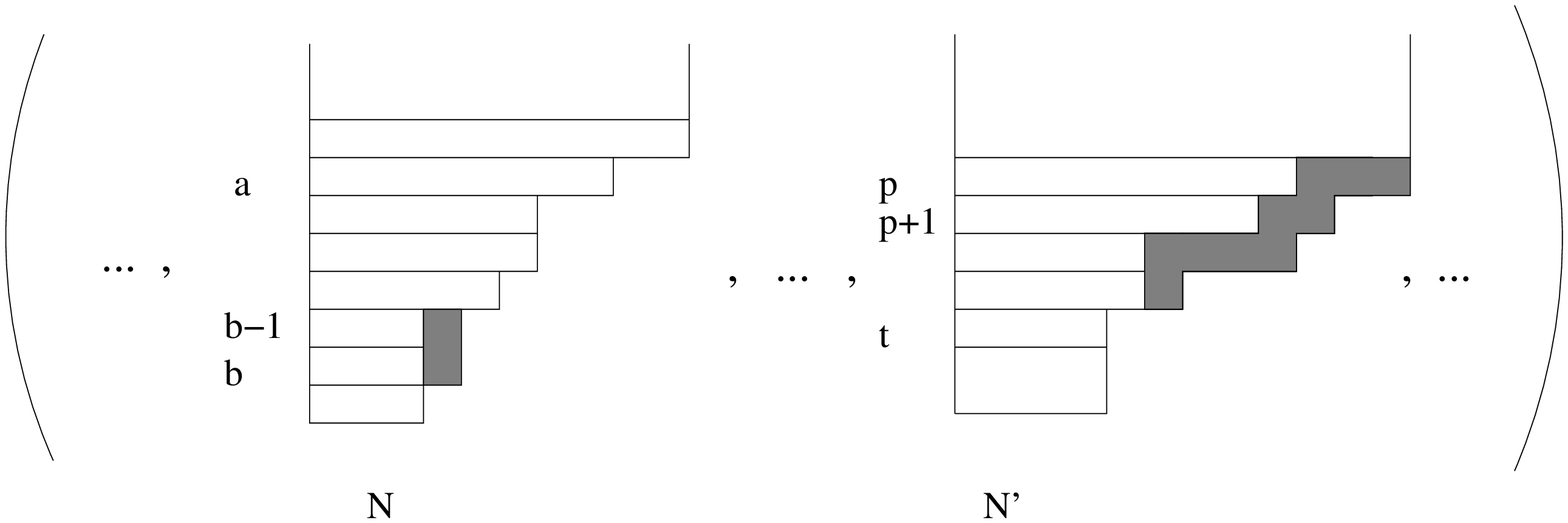}
\end{center}
 As we have $\nu_{t}^{(b')}-1+y_t=\nu_{t-1}^{(b')}$, note that:
$$\lambda[1]_{j_1}^{(b)}-j_1+s_{b}-\frac{(b-1)e}{l}>\nu_{t-1}^{(b')}-(t-1)+s_{b'}-\frac{(b'-1)e}{l}.\ \ \ \ \ \ (**)$$
We now define a $l$-composition $\underline{\lambda}[2]$ as follows. Put :
$$x_{j_1}':=x_{j_1}+x_{j_1+1}+...+x_s-y_t.$$
\begin{itemize}
\item If $y_{t-1}\leq x_{j_1}'$,  we put $j_2:=j_1$ and we define:
$$\underline{\lambda}[2]_i^{(k)}=\left\{ \begin{array}{ll}
{\lambda}[1]_i^{(k)}&\textrm {if }(i,k)\neq (j_1,b) \textrm{ and }(i,k)\neq (t-1,b'),\\
{\lambda}[1]_i^{(k)}+y_{t-1}&\textrm {if }(i,k)= (j_1,b),\\
{\lambda}[1]_i^{(k)}-y_{t-1}&\textrm {if }(i,k)= (t-1,b').\\
\end{array}\right.$$
By $(**)$, we have:
$${\lambda}[1]_{j_1}^{(b)}-j_1+s_b-\frac{(b-1)e}{l}>{\lambda}[1]_{t-1}^{(b')}-y_{t-1}-(t-1)+s_b-\frac{(b'-1)e}{l}.$$
By Proposition \ref{combi}, we get:
$$\underline{\lambda}[2]\prec\underline{\lambda}[1].$$

\item If  $y_{t-1}> x_{j_1}'$, there exists $j_2\in{\{r,r+1,...,j_1-1\}}$ such that:
$$x_{j_2}+x_{j_2+1}+...+x_{j_1-1}+x_{j_1}'\geq y_{t-1}>x_{j_2+1}+...+x_{j_1-1}+x_{j_1}'.$$
Then, we define:

$$\underline{\lambda}[2]_i^{(k)}=\left\{\begin{array}{ll}
{\lambda}[1]_i^{(k)}+x_{j_1}'&\textrm {if }(i,k)= (j_1,b),\\
{\lambda}[1]_i^{(k)}+x_{i}&\textrm {if }k=b\textrm{ and }\\
 & \textrm{   } j_1-1\geq i>j_2,\\
{\lambda}[1]_i^{(k)}-y_{t-1}&\textrm {if }(i,k)= (t-1,b'),\\
{\lambda}[1]_i^{(k)}+y_{t-1}-(x_{j_2+1}+...+x_{j_1-1}+x_{j_1}')&\textrm {if }(i,k)= (j_2,b),\\
{\lambda}[1]_i^{(k)}& \textrm {if otherwise.}
\end{array}\right.$$
By $(**)$, we have:
$${\lambda}[1]_{j_1}^{(b)}-i+s_{b}-\frac{(b-1)e}{l}>{\lambda}[1]_{t-1}^{(b')}-(t-1)-y_{t-1}+s_{b'}-\frac{(b'-1)e}{l}.$$
Moreover, we have  $\lambda[1]_{j_1}^{(b)}+x_{j_1}'=\lambda[1]_{j_1-1}^{(b)}+1$ and $\lambda[1]^{(b)}_i+x_i=\lambda[1]^{(b)}_{i-1}+1$ for all $i$ such that $j_1-1\geq i\geq r$. Thus, for all $i$ such that $j_1-1\geq i\geq j_2$, by $(**)$, we obtain:
$${\lambda}[1]_{i}^{(b)}-i+s_{b}-\frac{(b-1)e}{l}>{\lambda}[1]_{t-1}^{(b')}-(t-1)+x_{i+1}+...+x_{j_1-1}+x_{j_1}'-y_{t-1}+s_{b'}-\frac{(b'-1)e}{l}.$$
Hence, by Proposition \ref{combi}, we conclude that:
$$\underline{\lambda}[2]\prec\underline{\lambda}[1].$$
\end{itemize}
Keeping the above figure, $\underline{\lambda}[2]$ is  as follows.
\begin{center}
\leavevmode
\epsfxsize= 12cm
\epsffile{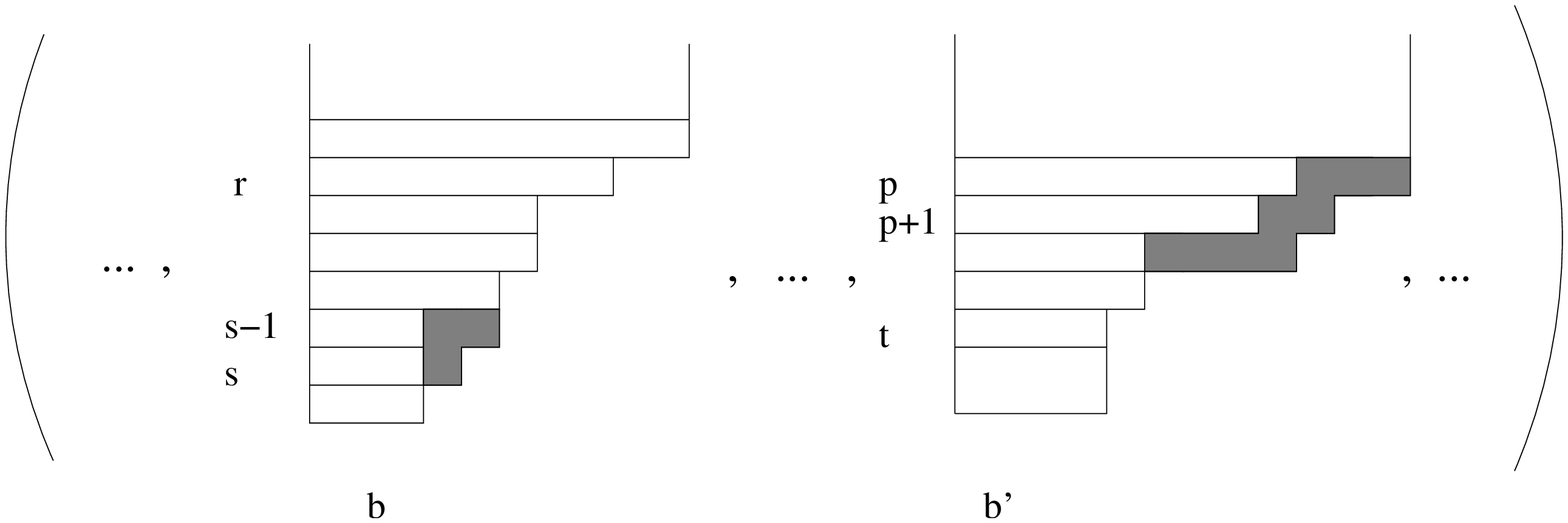}
\end{center}

Continuing in this way, we obtain the $l$-composition  $\underline{\lambda}[p-t+1]=\underline{\lambda}$:
\begin{center}
\leavevmode
\epsfxsize= 12cm
\epsffile{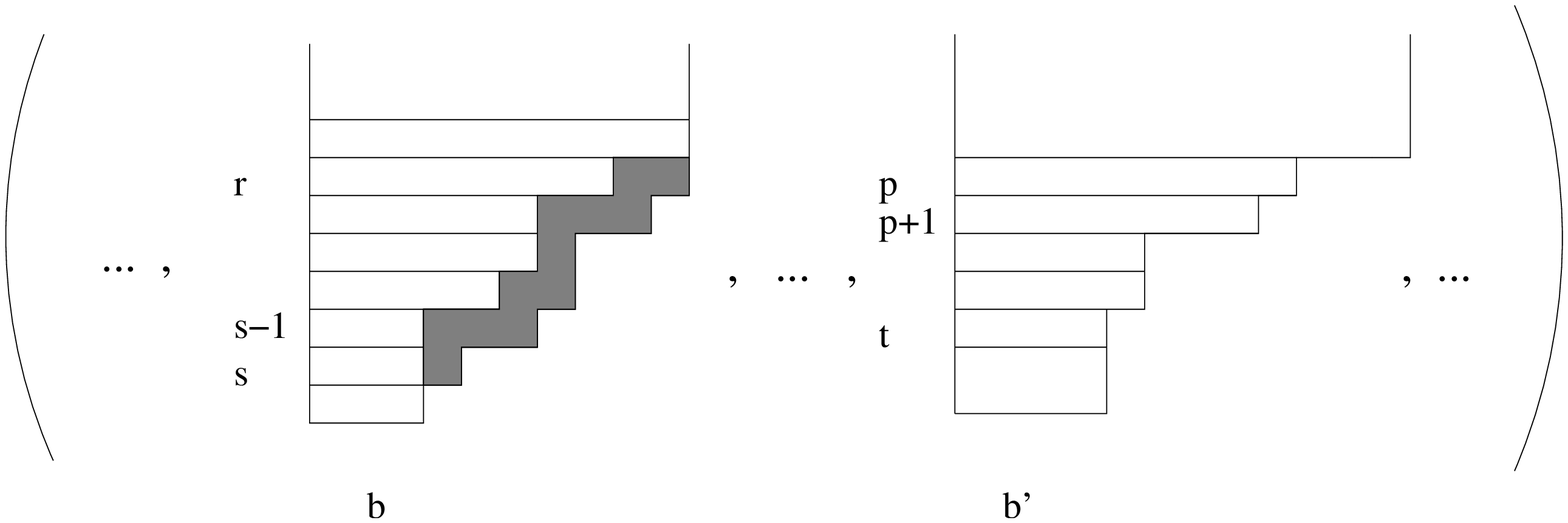}
\end{center}
We have :
$$\underline{\lambda}\prec\underline{\lambda}'.$$
Thus, we conclude that:
$$a_{{\bf s}_l}(\underline{\lambda})< a_{{\bf s}_l}(\underline{\lambda}').$$
Now  all the ordered monomials $u_{\underline{p}}$ which appear in $\overline{u_{\underline{k}}}$ are obtained recursively by using the relations $(R1)-(R4)$. Hence, by induction,  there exists $\beta$ such that:
$$\overline{u_{\underline{k}}}= \beta u_{\underline{k}}+\textrm{ sum of ordered monomials }u_{\underline{r}}\textrm{ with }a(u_{\underline{r}})>a(u_{\underline{k}}), $$
By \cite[Remark 3.24]{U}, we have $\beta=1$. This concludes the proof of the proposition.

\subsection{Consequences}

The result of the previous section leads to the following theorem:
\begin{thm}\label{mainm}
Let  ${\bf s}_l\in{\mathbb{Z}^n}$ and let $M_{{\bf s}_l}:=U_q' (\widehat{\mathfrak{sl}_e}).|\underline{\emptyset},{\bf s}_l\rangle\subset \mathcal{F}_{e,{\bf s}_l}$. Let
$$\{G(\underline{\mu},{\bf s}_l)\ |\ \ \underline{\mu}\in{\Lambda^{n}_{e;{\bf s}_l}   },\ n\in{\mathbb{N}}\}$$
be the canonical basis elements. Then, for all  $n\in{\mathbb{N}}$ and $\underline{\mu}\in{\Lambda^{n}_{e;{\bf s}_l}}$, we have:
$$G(\underline{\mu},{\bf s}_l)=|\underline{\mu},{\bf s}_l\rangle+\sum_{\underline{\lambda}\in{\Pi_l^n}, a_{{\bf s}_l}(\underline{\lambda})>a_{{\bf s}_l}(\underline{\mu})} d_{\underline{\lambda},\underline{\mu}}(q)|\underline{\lambda},{\bf s}_l\rangle,$$
where $d_{\underline{\lambda},\underline{\mu}}(q)\in q\mathbb{C}[q]$.
\end{thm}
\emph{Proof}.
 By Proposition \ref{main2}, we have for all $\underline{\lambda}\in{\Pi_l^n}$:
$$\overline{|\underline{\lambda},{\bf s}_l\rangle}=|\underline{\lambda},{\bf s}_l\rangle+\textrm{ sum of  } |\underline{\mu},{\bf s}_l\rangle \textrm{ with }
 a_{{\bf s}_l}(\underline{\mu})>a_{{\bf s}_l}(\underline{\lambda}).$$
Assume that  $\underline{\lambda}\in{\Lambda^{n}_{e;{\bf s}_l}}$ and that  $\underline{\mu}$ is one of the minimal $l$-partitions with respect to $a_{{\bf s}_l}$ which appears in $G(\underline{\lambda},{\bf s}_l)$. Then the characterization of the canonical basis in Theorem \ref{cb} immediately implies that   $\underline{\lambda}=\underline{\mu}$.\\
\rightline{$\Box$}

Now, we give the consequences  on the decomposition matrices of Ariki-Koike algebras. By Ariki's theorem, we obtain the following result. Note that this generalizes and gives a new proof of the main result of \cite{Jp} where the case $0\leq s_1\leq ...\leq s_l<e$ was considered (but the proof of \cite{Jp} gave an explicit construction of the canonical basis elements). Note also that it provides an interpretation of the parametrization by the Kleshchev $l$-partitions by using $a$-values.

Let ${\bf s}_l=(s_1,...,s_l)\in{\mathbb{Z}^l}$. Let $\mathcal{H}_{\mathbb{C},n}$ be the Ariki-Koike algebra over $\mathbb{C}$ with the following choice of parameters:
\begin{eqnarray*}
&&x_j=\eta_e^{s_{j}}\ \textrm{for}\ j=1,...,l,\\
&&x=\eta_{e},
\end{eqnarray*}
 
Let  $\mathcal{H}_{\mathbb{C}(y),n}^{{\bf s}_l}$ be the Ariki-Koike algebra over $\mathbb{C}(y)$   with the following parameters:
\begin{eqnarray*}
&&u_j=y^{lm^{(j)}}\eta_l^{j-1}\ \textrm{for}\ j=1,...,l,\\
&&v=y^l,
\end{eqnarray*}
where $\displaystyle m^{(j)}=s_j-\frac{(j-1)e}{l}+\alpha e$ for
$j=1,...,l$ and where $\alpha$ is a positive integer such that
$m^{(j)}>0$ for  $j=1,...,l$. Then the specialisation $\theta:\mathbb{C}[y,y^{-1}]\to \mathbb{C}$ which sends $y$ to $\displaystyle \textrm{exp}(\frac{2i\pi}{le})$ induces a decomposition map:
$$d:R_0 (\mathcal{H}^{{\bf s}_l}_{\mathbb{C}(y),n})\to R_0 (\mathcal{H}_{\mathbb{C},n}).$$
\begin{thm}\label{maing}  For each $M\in{\Irr{\mathcal{H}_{\mathbb{C},n}}}$, there exists a   unique simple $\mathcal{H}_{\mathbb{C}(y),n}^{{\bf s}_l}$-module $V_M$  such that the following two conditions hold:
\begin{itemize}
\item $d_{V_M,M}= 1$,
\item if there exists $W\in{\mathcal{H}_{\mathbb{C}(y),n}^{{\bf s}_l}}$ such that $d_{W,M}\neq 0$ then $a(W)>a(V_M)$.
\end{itemize}
  Moreover the function which sends $M$ to $V_M$ is injective. As a consequence the associated decomposition matrix is unitriangular and the following set is in natural bijection with $\Irr{\mathcal{H}_{\mathbb{C},n}}$:
$$\mathcal{B}_{{\bf s}_l}=\{V_M\ |\ M\in{\Irr{\mathcal{H}_{\mathbb{C},n}}}\}.$$
Finally, we have:
$$\mathcal{B}_{{\bf s}_l}=\{S^{\underline{\lambda}}_{\mathbb{C}(y)}\ |\ \underline{\lambda}\in{\Lambda^{n}_{e;{\bf s}_l}} \}.$$
\end{thm}
\begin{ex}\label{exfin}
The decomposition matrices of Ariki-Koike
algebras can be computed using the algorithm described in
\cite{Ja}. We put $l=2$,
$x=\eta_4$, $x_1=\eta_4^0$ and $x_2=\eta_4^1$. Then Theorem \ref{maing} shows
how we can order the rows of the associated decomposition matrix to
have a unitriangular shape with $1$ along the diagonal. To do this, we
must choose $(s_1,s_2)$ such that $s_1\equiv 0 (\textrm{mod }4)$ and
$s_1\equiv 1 (\textrm{mod }4)$ and order the rows with respect to the
$a$-values.  These rows are indexed by the $2$-partitions of rank $4$,
we also give the associated $a$-values. We obtain the explicit
bijection between  $\Lambda^{4}_{4,(s_1,s_2)}$ and $\Irr{\mathcal{H}_{\mathbb{C},n}}$.
\begin{itemize}
\item When $s_1=0$ and $s_2=1$, the simple $\mathcal{H}_{\mathbb{C},n}$-modules may be labeled by $\Lambda^{4}_{4,(0,1)}$.

$$\Lambda^4_{4,(0,1)}=\{(\emptyset,(4));((1),(2,1));((1,1),(1,1));((1),(3));((1,1),(2));$$
$$((2),(1,1));((2),(2));((2,1),(1));((2,1,1),\emptyset);((2,2),\emptyset);((3),(1));$$
$$((3,1),\emptyset);((4),\emptyset)\};$$
$$
\begin{array}{cc}
      ((4),\emptyset) & 0 \\
     ((3),(1))  & 1\\
     (\emptyset,(4)) & 1\\
     ((3,1),\emptyset)  &1\\
     ((2),(2))  &2\\
     ((2,2),\emptyset)   &2\\
     ((1),(3))  &2\\
   ((2,1),(1))  &3\\
    ((2,1,1),\emptyset)  &4\\
    ((2),(1,1))  &4\\
    ((1,1),(2))  &4\\
    ((1),(2,1))  &5\\
   ((1,1),(1,1)) &6 \\
     (\emptyset,(3,1))  &4\\
   ((1,1,1),(1))  &6\\
     (\emptyset,(2,2))   &6\\
   ((1,1,1,1),\emptyset)  &9\\
    (\emptyset,(2,1,1))  &9\\
   ((1),(1,1,1))  &9\\
  (\emptyset,(1,1,1,1))  &16
\end{array}
\left(
\begin{array}{cccccccccccccc}
  &1&.&.&.&.&.&.&.&.&.&.&.&.\\
  &.&1&.&.&.&.&.&.&.&.&.&.&.\\
  &.&.&1&.&.&.&.&.&.&.&.&.&.\\
  &1&.&.&1&.&.&.&.&.&.&.&.&.\\
  &.&1&.&.&1&.&.&.&.&.&.&.&.\\
  &.&.&.&.&.&1&.&.&.&.&.&.&.\\
 &1&.&1&.&.&.&1&.&.&.&.&.&.\\
 &.&.&.&.&.&.&.&1&.&.&.&.&.\\
  &.&.&.&1&.&.&.&.&1&.&.&.&.\\
  &.&.&.&.&.&.&.&.&.&1&.&.&.\\
  &1&.&.&1&.&.&1&.&.&.&1&.&.\\
  &.&.&.&.&.&.&.&.&.&.&.&1&.\\
  &.&.&.&.&.&1&.&.&.&.&.&.&1\\
  &.&.&1&.&.&.&1&.&.&.&.&.&.\\
  &.&.&.&1&.&.&.&.&1&.&1&.&.\\
  &.&.&.&.&1&.&.&.&.&.&.&.&.\\
  &.&.&.&.&.&.&.&.&1&.&.&.&.\\
  &.&.&.&.&.&.&1&.&.&.&1&.&.\\
  &.&.&.&.&.&.&.&.&.&.&.&.&1\\
  &.&.&.&.&.&.&.&.&.&.&1&.&.\\
\end{array}\right)$$

\item When $s_1=4$ and $s_2=1$, the simple $\mathcal{H}_{\mathbb{C},n}$-modules may be labeled by $\Lambda^{4}_{4,(4,1)}$.

$$\Lambda^4_{4,(4,1)}=\{((1,1,1),(1));((1),(2,1);((1,1),(1,1));((1),(3));((1,1),(2));$$
$$((2),(1,1));((2),(2));((2,1),(1));((2,1,1),\emptyset);((2,2),\emptyset);((3),(1));$$ $$((3,1),\emptyset);((4),\emptyset)\};$$

$$
\begin{array}{cc}
((4),\emptyset)& 0\\
((3,1),\emptyset)& 1\\
((2,2),\emptyset)& 2\\
((2,1,1),\emptyset)& 3\\
((3),(1))& 5\\
((2,1),(1))& 6\\
((1,1,1),(1))& 8\\
((2),(2))& 9\\
((1,1),(2))& 10\\
((2),(1,1))& 12\\
((1),(3))& 12\\
((1,1),(1,1))& 13\\
((1),(2,1))& 16\\
((1,1,1,1),\emptyset)& 6\\
(\emptyset,(4))& 14\\
(\emptyset,(3,1))& 19\\
(\emptyset,(2,2))& 22\\
(\emptyset,(2,1,1))& 25\\
((1),(1,1,1))& 21\\
(\emptyset,(1,1,1,1))& 32\\
\end{array}
\left(\begin{array}{cccccccccccccc}
1 & . & .&.&.&.&.&.&.&.&.&.&.\\
1&1&.&.&.&.&.&.&.&.&.&.&.\\
.&.&1&.&.&.&.&.&.&.&.&.&.\\
.&1&.&1&&.&.&.&.&.&.&.&.\\
.&.&.&.&1&.&.&.&.&.&.&.&.\\
.&.&.&.&.&1&.&.&.&.&.&.&.\\
.&1&.&1&.&.&1&.&.&.&.&.&.\\
.&.&.&.&1&.&.&1&.&.&.&.&.\\
1&1&.&.&.&.&1&.&1&.&.&.&.\\
.&.&.&.&.&.&.&.&.&1&.&.&.\\
1&.&.&.&.&.&.&.&1&.&1&.&.\\
.&.&1&.&.&.&.&.&.&.&.&1&.\\
.&.&.&.&.&.&.&.&.&.&.&.&1\\
.&.&.&1&.&.&.&.&.&.&.&.&.\\
.&.&.&.&.&.&.&.&.&.&1&.&.\\
.&.&.&.&.&.&.&.&1&.&1&.&.\\
.&.&.&.&.&.&.&1&.&.&.&.&.\\
.&.&.&.&.&.&1&.&1&.&.&.&.\\
.&.&.&.&.&.&.&.&.&.&.&1&.\\
.&.&.&.&.&.&1&.&.&.&.&.&.
\end{array}\right)
$$
\item When $s_1=0$ and $s_2=5$, the simple $\mathcal{H}_{\mathbb{C},n}$-modules may be labeled by $\Lambda^{4}_{4,(0,5)}$.

$$\Lambda^4_{4,(0,5)}=\{(\emptyset,(2,1,1));(\emptyset,(2,2)),((1),(1,1,1));((1,1),(1,1));((1,1),(2));$$ $$((2),(1,1));((1),(2,1));((1,1,1),(1));((1,1),(2));((1),(3));(\emptyset,(3,1));$$ $$(\emptyset,(4));((2,1),(1)),((2),(2))\}.$$

$$
\begin{array}{cc}
(\emptyset,(4))& 0\\
(\emptyset,(3,1))& 1\\
(\emptyset,(2,2))& 2\\
((1),(3))& 3\\
(\emptyset,(2,1,1))& 3\\
((1),(2,1))& 4\\
((2),(2))& 5\\
((1),(1,1,1))& 6\\
((2),(1,1))& 6\\
((1,1),(2))& 8\\
((1,1),(1,1))& 9\\
((2,1),(1))& 10\\
((1,1,1),(1))& 15\\
(\emptyset,(1,1,1,1))& 6\\
((4),\emptyset)& 6\\
((3),(1))& 6\\
((3,1),\emptyset)& 11\\
((2,2),\emptyset)& 14\\
((2,1,1),\emptyset)& 17\\
((1,1,1,1),\emptyset)& 24\\
\end{array}
\left(\begin{array}{cccccccccccccc}

 1&. &. &.&. &. &. &. &. &. &. &. &.\\
 1&1&.&.&.&.&.&.&.&.&.&.&.\\
 .&.&1&.&.&.&.&.&.&.&.&.&.\\
1&1&.&1&.&.&.&.&.&.&.&.&.\\
.&1&.&.&1&.&.&.&.&.&.&.&.\\
.&.&.&.&.&1&.&.&.&.&.&.&.\\
.&.&1&.&.&.&1&.&.&.&.&.&.\\
.&.&.&.&.&.&.&1&.&.&.&.&.\\
.&.&.&.&.&.&.&.&1&.&.&.&.\\
.&1&.&1&1&.&.&.&.&1&.&.&.\\
.&.&.&.&.&.&.&1&.&.&1&.&.\\
.&.&.&.&.&.&.&.&.&.&.&1&.\\
.&.&.&.&1&.&.&.&.&1&.&.&1\\
.&.&.&1&.&.&.&.&.&.&.&.&.\\
.&.&.&.&.&.&1&.&.&.&.&.&.\\
.&.&.&.&1&.&.&.&.&.&.&.&.\\
.&.&.&1&.&.&.&.&.&1&.&.&.\\
.&.&.&.&.&.&.&.&.&.&1&.&.\\
.&.&.&.&.&.&.&.&.&1&.&.&1\\
.&.&.&.&.&.&.&.&.&.&.&.&1
\end{array}\right)
$$

\end{itemize}
\end{ex}

\bibliographystyle{amsalpha}

\begin{thebibliography}{A}

\bibitem[As]{As}
S. {Ariki},
\textit{ On the semi-simplicity of the Hecke algebra of $(\mathbf{Z}/r\mathbf{Z})\wr{{S}_n}$}.
J. Algebra \textbf{169}, no. 1 (1994): 216--225.



\bibitem[Ad]{Ad}
S. Ariki,
\textit{On the decomposition numbers of the Hecke algebra of $G(m,1,n)$}.
J. Math. Kyoto Univ., \textbf{36} (1996): 789--808.



\bibitem[Ac]{Ac}
S. Ariki,
\textit{On the classification of simple modules for Cyclotomic {H}ecke algebras of type ${G}(m,1,n)$ and Kleshchev multipartitions}.
Osaka J.Math., \textbf{38} (2001): 827-837.

\bibitem[Ab]{Ab} S. Ariki,
\textit{Representations of quantum algebras and combinatorics of Young tableaux}.
Univ. Lecture Series,
\textbf{26} (2002), AMS.


\bibitem[AK]{AK}
S. Ariki, K. Koike,
\textit{A Hecke algebra of $(\mathbb{Z}/r\mathbb{Z})\wr{\mathfrak{S}_n}$
and construction of irreducible representations}.
Adv. Math., \textbf{106} (1994): 216--243.

\bibitem[BM]{BM}
M. {Brou\'e}, G. {Malle},
\textit{Zyklotomische Heckealgebren}.
Ast\'erisque, \textbf{212} (1993): 119--189.


\bibitem[DJM]{DJM} R. Dipper, G. James, A. Mathas,
\textit{ Cyclotomic $q$-Schur algebras}.
Math. Z., \textbf{229} no. 3 (1998): 385--416.



\bibitem[DM]{DM} R. Dipper, A. Mathas, \textit{Morita equivalences of Ariki-Koike algebras}, Math.Z., \textbf{240} no. 3 (2003): 579--610.


\bibitem[FL]{FL}
 O. Foda, B. Leclerc, M. Okado, J-Y Thibon, T. Welsh,
\textit{  Branching functions of $A\sp {(1)}\sb {n-1}$ and Jantzen-Seitz problem for Ariki-Koike   algebras}.
 Adv. Math., \textbf{141} no. 2 (1999), 322--365.

\bibitem[GL]{GL} J. Graham, G. Lehrer,
\textit{Cellular algebras}.
Invent. Math., \textbf{123} (1996): 1--34.

\bibitem[Gm]{Gm}
M. {Geck},
\textit{Representations of Hecke algebras at roots of unity}.
 S\'eminaire Bourbaki. Vol. 1997/98. Ast\'erisque \textbf{252}  (1993): 33--55.

\bibitem [G]{G} M. Geck, \textit{Kazhdan-Lusztig cells and decompositions numbers}, Representation theory, \textbf{2} (1998): 264--277.

\bibitem [Gs]{Gs} M. Geck, \textit{Modular Representations of Hecke
  algebras}, EPFL Press (to appear).


\bibitem[GIM]{GIM}
M. {Geck}, L. {Iancu}, G. {Malle},
\textit{Weights of Markov traces and generic degrees}.
Indag. Math., \textbf{11} (2000): 379--397.


\bibitem[GR]{GR}
 M. Geck, R. Rouquier,
\textit{Filtrations on projective modules for {I}wahori-{H}ecke algebras}.
Modular representation theory of finite groups (Charlottesville, VA, 1998), de Gruyter, Berlin (2001): 211--221.

\bibitem[Jp]{Jp} N. Jacon,
\textit{On the parametrization of the simple modules for Ariki-Koike algebras},  J. Math. Kyoto Univ, {\bf 44} no. 4 (2004): 729--767.

\bibitem[Ja]{Ja} N. Jacon,
\textit{An algorithm for the computation of the decomposition matrices
  for Ariki-Koike algebras},  J. Algebra (Comp. Algebra section), {\bf
  292} (2005), 100--109.

\bibitem[JM]{JM}
M. Jimbo, K. {Misra}, T. {Miwa}, M. {Okado},
\textit{Combinatorics of representations of $\mathcal{U}_q(\widehat{sl}(n))$ at $q=0$}. Commun. Math. Phys., 136 (1991): 543--566.

\bibitem[LT]{LT}
B. Leclerc, J-Y Thibon,
\textit{Canonical bases of $q$-deformed Fock spaces}. Int. Math. Res. Notices, {\bf 9}  (1996): 447--456.


\bibitem[Ms]{Ms}
A. {Mathas}
\textit{Simple modules of Ariki-Koike algebras}, in Group representations: cohomology, group actions and topology, Proc. Sym. Pure Math., \textbf{63} (1998): 383--396.

\bibitem[Ma]{Ma}
A. {Mathas}
\textit{The representation theory of the Ariki-Koike and cyclotomic q-Schur algebras}, Representation theory of algebraic groups and quantum groups, Adv. Studies Pure Math. (2004): 261--320.




\bibitem[TU]{TU} H. Takemura,  D. Uglov,
\textit{Representations of the quantum toroidal algebra on highest weight modules of the quantum affine algebra of type ${\mathfrak{gl}}_N$.}
Publ. Res. Inst. Math. Sci. \textbf{35}, No.3  (1999): 407--450.

\bibitem[U]{U} D. Uglov,
\textit{Canonical bases of higher-level $q$-deformed Fock spaces and Kazhdan-Lusztig polynomials},
 Kashiwara, Masaki (ed.) et al., Boston: Birkh\"auser. Prog. Math. \textbf{191}  (2000): 249--299.

\bibitem[Us]{Us} D. Uglov,
\textit{Canonical base of higher-level q-deformed Fock spaces} (short version of \cite{U}) {\bf math.QA/9901032}.


\bibitem[Y]{Y}
X. Yvonne,
\textit{A conjecture for $q$-decomposition matrices of cyclotomic
  $v$-Schur algebras}, preprint  {\bf math.RT/0505379}.
\bibitem[Yt]{Yt}
X. Yvonne,
\textit{Bases canoniques d'espaces de Fock en niveau sup\'erieur}, PhD
thesis, Universit\'e de Caen (2005).


\end{thebibliography}

\end{document}